\documentclass[a4paper,10pt,final,notitlepage]{article}
\usepackage[latin1]{inputenc}
\usepackage{bbm}
\usepackage{amsmath}
\usepackage{amsthm}
\usepackage{amscd}
\usepackage{amsfonts}
\usepackage{amssymb}

\newtheorem{Theorem}{Theorem}[section]
\newtheorem{Definition}[Theorem]{Definition}
\newtheorem{Proposition}[Theorem]{Proposition}
\newtheorem{Lemma}[Theorem]{Lemma}

\newtheorem{Remark}[Theorem]{Remark}

\newtheorem{Notation}[Theorem]{Notation}

\newcommand{\ud}{\mathrm{d}}
\newcommand{\be}{\begin{equation}}
\newcommand{\ee}{\end{equation}}
\newcommand{\bd}{\begin{displaymath}}
\newcommand{\ed}{\end{displaymath}}
\newcommand{\nd}{\stackrel{def}{=}}
\newcommand{\s}{\hspace{.05in}}

\newcommand{\idu}{\mathbbmss{1}}
\newcommand{\id}{\mathrm{I}}
\newcommand{\lla}{\left\langle}
\newcommand{\rra}{\right\rangle}

\def\qqed{\hfill\hbox{\hskip 6pt\vrule width6pt height7pt
depth1pt  \hskip1pt}\bigskip}

\author{G. Fabbri \footnote{Dipartimento di matematica, Universit\'a "La Sapienza", Roma, Italy {\tt fabbri@mat.uniroma1.it}}}
\title{A viscosity solution approach to the infinite dimensional HJB equation related to boundary control problem in transport equation}
\begin{document}

\maketitle

\begin{abstract}
The paper concerns with the infinite dimensional Hamilton-Jacobi-Bellman equation related to optimal control problem regulated by a linear transport equation with boundary control. 
A suitable viscosity solution approach is needed in view of the presence of the unbounded control-related term in the state equation in Hilbert setting. An existence-and-uniqueness result is obtained.


\end{abstract}

\section*{Introduction}
We study the Hamilton-Jacobi-Bellman equation (from now HJB equation) related to the infinite dimensional formulation of an optimal control problem whose state equation is a PDE of transport type.

We consider the PDE
\be
\label{introPDE}
\left\{ \begin{array}{ll}
\frac{\partial}{\partial s} x(s,r) +\beta\frac{\partial}{\partial r} x(s,r) = -\mu x(s,r) + \alpha(s,r) \;\;\; (s,r)\in (0,+\infty)\times(0,\bar s) \\
x(s,0)= a(s) \s\;\;\;\; if \;\s s > 0\\ 
x(0,r)= x^0(r) \s\;\;\;\; if \s\; r\in[0,\bar{s}]
\end{array} \right.
\ee
where $\bar s,\beta$ are positive constants, $\mu\in\mathbb{R}$, the initial data $x^0$ is in $L^2(0,\bar s)$ and we consider two controls: a boundary control $a$ is in $L^2_{loc} ([0, +\infty);\mathbb{R})$ and a distributed control $\alpha \in L^2_{loc}([0,+\infty)\times [0,\bar s]; \mathbb{R})$. More precisely we will consider controls such that $\alpha(\cdot)\in\mathcal{E}$ and $a(\cdot)\in\mathcal{A}$ where $\mathcal{E}$ and $\mathcal{A}$ will be defined in Section \ref{sezionenotazione}. We write ``$-\mu x$'' instead of ``$\mu x$'' because it is the standard way to write the equation in the economic literature where $-\mu$ has the meaning of a depreciation factor (and only the case $\mu\geq 0$ is used). Here we consider a generic $\mu \in \mathbb{R}$. 

Using the approach and the references described in Section \ref{sezionenotazione}, the above equation can be written as an ordinary differential equation in the Hilbert space $\mathcal{H}=L^2(0,\bar s)$ as follows
\be
\label{introeqstatoinDAstarprimo}
\left\{ \begin{array}{ll}
\frac{\ud}{\ud s} x(s) = A x(s) - \mu x(s) + \alpha(s) + \beta \delta_0 a(s)\\
x(0)=x^0
\end{array} \right.
\ee
where $A$ is the generator of a suitable $C_0$ semigroup and $\delta_0$ is the Dirac delta in $0$. Such an unbounded contribution in the Hilbert formulation comes from the presence in the PDE of a boundary control (see \cite{BDDM}).
%and the boundary control that appears in the PDE, as usual (see \cite{BDDM}), gives a non-bounded contribution in the Hilbert formulation, emphasized in our expression by the Dirac's symbol $\delta_0$.
We consider the problem of minimizing the cost functional
\be
\label{introtargetfunctional}
J(x,\alpha(\cdot),a(\cdot))=\int_0^\infty e^{-\rho s} L(x(s),\alpha(s),a(s)) \ud s
\ee
where $\rho>0$ and $L$ is globally bounded and satisfies some Lipschitz-type condition, as better described in Section \ref{sezionenotazione}. The HJB equation related to the control problem with state equation (\ref{introeqstatoinDAstarprimo}) and target functional (\ref{introtargetfunctional}) is
\begin{multline}
\label{introHJB}
\rho u(x) - \lla \nabla u(x) , Ax \rra - \lla \nabla u, -\mu x \rra_{L^2(0,\bar s)} -\\
- \inf_{(\alpha,a)\in\Sigma\times\Gamma} \left( \lla \beta \delta_0 (\nabla u),a \rra_\mathbb{R} + \lla \nabla u, \alpha \rra_{L^2(0,\bar s)} +  L(x,\alpha,a) \right)  = 0.
\end{multline}
The sets $\Gamma$ and $\Sigma$ will be introduced in Section \ref{sezionenotazione}. If we define the value function of the control problem as
\bd
V(x) \nd \inf_{(\alpha(\cdot),a(\cdot)) \in \mathcal{E}\times \mathcal{A}} J(x,\alpha(\cdot),a(\cdot)),
\ed
we wish to prove that it is the only solution, in a suitable sense, of the HJB equation.

We use the viscosity approach. Our main problem is to write a suitable definition of viscosity solution, so that an existence and uniqueness theorem can be derived for such a solution. The main difficulties we encounter, with respect to the existing literature, is in dealing with the boundary term and the non-analyticity of the semigroup. We substantially follow the original idea of Crandall and Lions (\cite{CL4} and \cite{CL5}) - with some changes, as the reader will rate in Definition \ref{defsubsol} and Definition \ref{defsupersol} - of writing test functions as the sum of a \textquotedblleft good part\textquotedblright\phantom{ }, a regular function with differential in $D(A^*)$ and a \textquotedblleft bad part\textquotedblright\phantom{ } represented by some radial function. The main problems arise from the evaluation of the boundary term on the radial part.

In order to write a working definition in our case some further requirements are needed, like a $C^2$ regularity of the test functions, the presence of a \textquotedblleft remainder term\textquotedblright\phantom{ } in the definition of sub/super solution and the $B$-Lipschitz continuity (see Definition \ref{defblipschitz}) of the solution. This last feature guarantees that the maxima and the minima in the definition of sub/super solution remain in $D(A^*)$ (see Proposition \ref{prop0}). Some other comments on the definition of solution (Definition \ref{defsubsol} and \ref{defsupersol}) need some technical details and can be found in Remark \ref{remarkdefinizione}.

The used technique cannot be easily extended to treat a general non-linear problem because we use the explicit form of the PDE that we give in (\ref{eq5bis}). Moreover the (linear) case in which $\mu$ is not a constant but a function of $L^{\infty}(0,\bar s)$ (and so we have $\mu(r)$ in equation (\ref{introPDE})) presents some difficulties in the proof of the uniqueness theorem and so it is not an easy generalization of the present work. See Remark \ref{remarkmunondipendedar} for details.

\paragraph*{A brief summary of the literature}
Hamilton-Jacobi equations in infinite dimensions, especially when arising from optimal control problems in Hilbert spaces, was first studied by Barbu and Da Prato (\cite{barbudaprato1}, \cite{barbudaprato2}) with strong solutions approach. The viscosity method, introduced in the study of finite dimensional HJ equations in \cite{CLfinitedim} was generalized by the same authors in a series of works; the most important for our approach are \cite{CL4} and \cite{CL5}. Moreover new variants of the notion of viscosity solution for HJB equations in Hilbert space was given in \cite{ishii}, \cite{tataru1}, \cite{tataru2}, \cite{tataru3} and \cite{CL6}.

The study of viscosity solution for HJB equations in Hilbert spaces arising from optimal control problem of systems modeled by PDE with boundary control term is more recent. In this research field there is not an organic and complete theory but some works on specific PDE that adapt the ideas and the techniques of viscosity solutions to particular problems using their own characteristics like we do in this work for the problem regulated by transport equation. For the first order HJB equations see \cite{CannarsaGozziSoner}, \cite{CannarsaTessitore} (see also \cite{CannarsaTessitore96,CannarsaTessitore96b}) where some classes of parabolic equations are treat,  \cite{GozziSritharanSwiech} in which the authors study the HJB equation related to a two-dimensional Navier-Stokes (see also \cite{Shimano02}). It must be noted that all these works  treat the case of $A$ analytic.

HJB equation like (\ref{introHJB}) was treated, only in the convex case, with strong solutions approach adapting the Barbu and Da Prato' arguments in \cite{Silviatesi} and \cite{faggian2}.

\paragraph*{A motivating economic problem}
Tranport equations are used to model a large variety of phenomena, from age-structured population models (see for instance \cite{Iannelli95, Anita01, Iannelli06}) to population economics (\cite{FeichtingerPrskwetzVeliov04}), from epidemiologic studies to socio-economic science and transport phenomena in phisics.

Problems such as (\ref{introPDE}) can be used to describe, in economics, capital accumulation processes where an heterogeneous capital is involved, and this is the reason why the study of infinite dimensional control problem is of growing interest in the economic fields. For instance in the vintage capital models $x(t,s)$ may be regarded as the stock of capital goods differentiated with respect the time $t$ and the vintage $s$. Heterogeneous capital, both in finite and infinite dimensional approach, is used to study depreciation and obsolescence of physical capital, geographical difference in growth, innovation and R\&D. 

Regarding problems modeled by a transport equation where an infinite dimensional setting is used we cite the following papers: \cite{BarucciGozzi98} and \cite{BarucciGozzi01} on optimal technology adoption in a vintage capital context (in the case of quadratic cost functional), \cite{hartl} on capital accumulation, \cite{BarucciGozzi99} on optimal advertising and \cite{faggian2} \cite{faggian3} that studies the convex functional case using a strong solutions approach. See also \cite{FaggianGozzi04}.

Moreover, we mention that the infinite dimensional approach may apply to problems such as issuance of public debt (see \cite{piccoli} for a description of the problem). In that problem a stochastic setting and simple state-control constraints appear, but hopefully the present work can be a first step in this direction.
% To understand, in a qualitative way, why an Hilbertian approach could be interesting in the public debt problem is enough to observe that at every time the new issuance of public bonds have to refinance expiring bonds, so the "state" of the system involve all the "history" of the issuances (and the history of the issuances is nothing else than an $L^2$ function).

\paragraph*{Plan of the paper}
The work is organized as follows: in the first section we remind some results on the state equation, we introduce some preliminary remarks on the main operators involved in the problem, we explain some notations, we define the HJB equation and we give the definition of solution. The second section regards some properties of the value function (in particular some regularity properties) that we will be used in the third section to prove that it is the only (viscosity) solution of the HJB equation.

\paragraph*{Acknowledgements} The author would like to thank Prof. Andrzej \'Swi\c ech for his hospitality, his great kindness, for many useful suggestions and stimulating conversations.

\section{Notation and preliminary results}
\label{sezionenotazione}
\subsection{State equation}
\label{sezioneequazione}
In this subsection we will see some properties of the state equation: we write it in three different (and equivalent) forms that point out different properties of the solution. We will use all the three forms in the following proofs.

We consider the PDE on $[0, +\infty) \times [0,\bar s]$ given by
\be
\label{eqstatoderivparziali}
\left\{ \begin{array}{ll}
\frac{\partial}{\partial s} x(s,r) +\beta\frac{\partial}{\partial r} x(s,r) = -\mu x(s,r) + \alpha(s,r) \;\;\; (s,r)\in (0,+\infty)\times(0,\bar s)\\
x(s,0)= a(s) \s\;\;\;\; if \s\; s > 0\\
x(0,r)= x^0(r) \s\;\;\;\; if \s\; r\in[0,\bar{s}]
\end{array} \right.
\ee
Given an initial datum $x^0 \in L^2((0,\bar s);\mathbb{R})$ (from now simply $L^2(0,\bar s)$), a boundary control $a(\cdot) \in L^{2}_{loc}([0, +\infty);\mathbb{R})$
%\footnote{We could take $a(\cdot) \in L^{2}_{loc}([0, +\infty);\mathbb{R})$ but in the sequel we will take only bounded control so we are interested only in the case $a(\cdot) \in L^{\infty}([0, +\infty);\mathbb{R})$} 
and a distributed control $\alpha(\cdot)\in L^{2}_{loc}([0,+\infty)\times[0,\bar s];\mathbb{R})$ the (\ref{eqstatoderivparziali}) has a unique solution in $L^{2}_{loc}([0,+\infty) \times [0, \bar s];\mathbb{R})$ given by
\be
\label{eq5bis}
x(s,r)=
\left\{ \begin{array}{ll}
e^{-\mu s}x^0(r-\beta s) + \int_0^s e^{-\mu \tau} \alpha(s-\tau, r -\beta \tau) \ud \tau & \; r\in[\beta s, \bar s]\\
e^{\frac{-\mu}{\beta} r}a(s-r/\beta) + \int_0^{r/\beta} e^{-\mu \tau} \alpha(s-\tau, r - \beta \tau) \ud \tau  & \; r\in [0,\beta s)
\end{array} \right.
\ee
In the following $x(s,r)$ will denote (\ref{eq5bis}).

We can rewrite such equation in a suitable Hilbert space setting. We take the Hilbert space $\mathcal{H} \nd L^2(0,\bar{s})$ and the $C_0$ semigroup $T(t)$ given by
\bd
T(s)f[r] \nd \left\{ \begin{array}{ll}
f(r-\beta s) & for \s r \in [\beta s,\bar{s}]\\
0 & for \s r \in [0,\beta s)
\end{array} \right.
\ed
The generator of $T(s)$ is the operator $A$ given by
\bd \left\{
\begin{array}{ll}
D(A)= \{ f \in H^1[0,\bar{s}] \s : \s f(0) = 0 \} \\
A(f)[r]=-\beta \frac{\ud}{\ud r} f (r)
\end{array} \right.
\ed
(see \cite{BarucciGozzi01} for a proof in the case $\beta=1$, the proof in our case can be obtined simply taking $s'=\beta s$).\\
\begin{Remark}
\label{notationquadre}
To avoid confusion if $x\in L^2(0,\bar s)$ we will use $[\cdot]$ to denote the pointwise evaluation, so $x[r]$ is the value of $x$ in $r\in[0,\bar s]$. On other hand $x(s)$ will denote the evolution of the solution of the state equation (in the Hilbert space) at time $s$ (as in (\ref{eq5tris})). That is, $x(s)$ is an element of $\mathcal{H}$ while $x[r]$ is a real number.\\
\end{Remark}

We want to write an infinite dimensional formulation of (\ref{eqstatoderivparziali}) but in $L^2(0,\bar s)$ it should appear like
\be
\label{eq5tris}
\left\{ \begin{array}{ll}
\frac{\ud}{\ud s} x(s) = A x(s) -\mu x(s) +\alpha(s) + \beta \delta_0 a(s)\\
x(0)=x^0
\end{array} \right.
\ee
where $\alpha(s)\in L^2(0, \bar s)$ is the function $r\mapsto \alpha(s,r)$. Such expression does not make sense in $L^2(0,\bar s)$ for the presence of the unbounded term $\beta\delta_0 a(s)$. We can anyway apply formally the variation of constants method to (\ref{eq5tris}) and obtain a mild form of (\ref{eq5tris}) that is continuous $\colon [0,+\infty) \to L^2(0,\bar s)$. This is what we do in the next definition. Note that we have written (\ref{eq5tris}) only to be more clear but we could go on in a more formal way without it.
\begin{Definition}
Given $x^0\in L^2(0,\bar{s})$, $a(\cdot) \in L^{2}_{loc}([0, +\infty);\mathbb{R})$ and \linebreak[4]
$\alpha(\cdot)\in L^{2}_{loc}([0,+\infty); L^2(0,\bar s))$ the function in $C([0,+\infty); L^2(0,\bar{s}))$ given by
\begin{multline}
\label{eqstatomildform}
x(s)= e^{-\mu s} T(s) x^0 - A \int_0^s e^{-\mu(s-\tau)} T(s-\tau) ( a(\tau) \nu) \ud \tau + \\
+\int_0^s e^{-\mu(s-\tau)} T(s-\tau) \alpha(\tau) \ud \tau
\end{multline}
where 
\bd
\begin{array}{ll}
\nu \colon [0,\bar{s}] \to \mathbb{R}\\
\nu \colon r\mapsto e^{-\frac{\mu}{\beta}r}
\end{array}
\ed
is called mild solution of (\ref{eq5tris}). 
\end{Definition}

\begin{Remark}
We could include the term $-\mu x$ in the generator of the semigroup $A$ taking a $\tilde A=A-\mu \idu$. In this case the equation (\ref{eq5tris}) would appear in the following equivalent form:
\begin{equation}
\label{eqstatoconatilde}
\frac{\ud}{\ud s} x(s) = \tilde A x(s) +\alpha(s) + \beta \delta_0 a(s)
\end{equation}
The problem of this approach is that often we will use, in the estimates, the dissipativity of the generator and $\tilde A$ is dissipative only if $\mu \geq 0$.
Nevertheless we wrote the mild form (\ref{eqstatomildform}) as if we consider the (of course equivalent!) state equation (\ref{eqstatoconatilde}), indeed in the definition we consider for example the convolution term given by
\[
\int_0^s e^{-\mu(s-\tau)} T(s-\tau) \alpha(\tau) \ud \tau
\]
and $e^{-\mu(s)} T(s)$ is exactly the semigroup generate by $\tilde A$. In the sequel (see (\ref{eqstatomildformdiversa})) we will use also another mild form of the state equation that is not explicit and it is the one obtined if we do not include the term $e^{-\mu s}$ in the semigroup. The two forms are equivalent.
\end{Remark}

\begin{Proposition}
Taken $x(s)$ the function $\colon \mathbb{R}^+ \to L^2(0,\bar s)$ given by (\ref{eqstatomildform}) and $x(s,r)$ the function $\colon \mathbb{R}^+\times [0,\bar s] \to \mathbb{R}$ given by (\ref{eq5bis}) we have $x(s)[r]=x(s,r)$.
\end{Proposition}
\begin{proof}
%--- questo non e' troppo vero, guardare Barucc e Gozzi 97, Iannelli e forse LasieckaTriggiani
See \cite{BarucciGozzi01}.
\end{proof}
Eventually we observe that (\ref{eq5tris}) can be rewritten in a precise way in a larger space where we have not problem with the $\beta\delta_0$ term. The mild solution will be the only solution (in a suitable sense) of the new differential equation. We need more notation to write it.

We consider the adjoint operator $A^*$. Its explicit expression is given by
\bd
\left\{
\begin{array}{ll}
D(A^*)\nd\{ f \in H^1(0,\bar{s}) \s :\s f(\bar{s}) =0 \} \\
A^*(f)[r]=\beta \frac{\ud}{\ud r} f (r)
\end{array}
\right.
\ed
On $D(A^*)$ we put the graph norm and the related Hilbert structure. We consider the inclusion
\bd
i\colon D(A^*) \hookrightarrow L^2(0,\bar{s})\\
\ed
and its continuous adjoint
\bd
i^*\colon L^2(0,\bar{s}) \to D(A^*)' 
\ed
where we have identified $L^2$ with its dual.

We can extend $A$ to a generator of a $C_0$ semigroup on $D(A^*)'$ (the domain of the extension will contain $L^2$) and we observe that the Dirac's measure $\delta_0 \in D(A^*)'$ (see \cite{Silviatesi} Proposition 4.5 page 60 for details).
\begin{Proposition}
Given $T>0$, $x^0 \in L^2(0,\bar{s})$, $a(\cdot)\in L^2(0,T)$, $\alpha(\cdot) \in L^2([0,T];L^2(0,\bar s))$, (\ref{eqstatomildform}) is the unique solution of 
\be
\label{eqstatoinDAstarprimo}
\left\{ \begin{array}{ll}
\frac{\ud}{\ud s} i^* x(s) = A x(s) -\mu x(s) + \alpha(s) + \beta\delta_0 a(s)\\
x(0)=x^0
\end{array} \right.
\ee
in $W^{1,2}(0,T; D(A^*)') \cap C(0,T, \mathcal{H})$. 
Moreover if $a(\cdot) \in W^{1,2}(0,T)$ then such solution will belong to $C^1(0,T; D(A^*)') \cap C(0,T; \mathcal{H})$.
\end{Proposition}
\begin{proof}
%--- questa e' proprio falsa. dovrebbe essere capitolo 3.2
See \cite{BDDM} Chapter 3.2 (in particular Theorem 3.1 page 173).
\end{proof}

\subsection{The definition of the operator $B$}
In this subsection we give the definition of the operator $B$ that will have a fundamental role. We could use an abstract approach, noting that $A$ and $A^*$ are both generator of $C_0$ semigroups of contractions and then both are negative (see \cite{DaPratoZabczyk92} page 424) and the set $\{ \lambda \in \mathbb{C} \s : \s \mathbb{R}e(\lambda) >0 \}$ is in the resolvent of both $A$ and $A^*$ (Hille-Yosida theorem, see \cite{LiYong} page 53). Anyway in this case we can also follow a more direct approach that allows to find the explicit form of the operator.

To note that that $A^*$ and $A$ are negative operators we take $\phi\in D(A^*)$ (so $\phi(\bar s)=0$)
\bd
\lla A^* \phi, \phi \rra = \int_0^{\bar s} \beta \phi'(r) \phi(r) \ud r = \frac{-\beta \phi(0)^2}{2}
\ed
and for $\phi\in D(A)$ (so $\phi(0)=0$)
\bd
\lla A \phi, \phi \rra = \int_0^{\bar s} -\beta \phi'(r) \phi(r) \ud r = \frac{-\beta \phi(\bar s)^2}{2}
\ed
So, given a $\lambda >0$, the operators $(A-\lambda \id)$ and $(A^* - \lambda \id)$ are strongly negative: $\lla (A- \lambda \id) x,x \rra \leq -\lambda |x|_{\mathcal{H}}^2$ for all $x\in D(A)$ and $\lla (A^*- \lambda \id) x,x \rra \leq -\lambda |x|_{\mathcal{H}}^2$ for all $x\in D(A^*)$.

We can also directly prove that
\bd
(A-\lambda \id)^{-1} \colon \mathcal{H} \to D(A)
\ed
is a continuous negative linear operator whose explicit expression is given by
\bd
(A-\lambda \id)^{-1}(\phi)[r] = \frac{1}{\beta} \left ( -e^{-\frac{\lambda}{\beta}r} \int_0^r e^{\frac{\lambda}{\beta}\tau} \phi(\tau) \ud \tau  \right )
\ed
%To check the continuity we observe that
%\begin{multline}
%|(A -\lambda \id)^{-1}(\phi)|^2_{\mathcal{H}} = \int_0^{\bar s} \left | e^{-\lambda r} \int_0^{r}  e^{\lambda \tau} \phi(\tau) \ud \tau  \right |^2 \ud r \leq\\
%\leq \int_0^{\bar s} e^{-\lambda r} \left ( \int_0^{r} e^{\lambda r} |\phi(\tau)| \ud \tau  \right )^2 \ud r\leq C_{\bar s} \int_0^{\bar s} \int_0^{\bar s} e^{2 \lambda \bar s} \phi^2(\tau) \ud \tau \ud r = \tilde C_{\bar s} |\phi|^2_{\mathcal{H}}
%\end{multline}
%where $C_{\bar s}$ and $\tilde C_{\bar s}$ are constants that depend only on $\bar s$ and $\lambda$. 
The continuity can be proven directly with not difficult estimates and the negativity can be proven directly using an integration by part argument. 

%Calculating explicitly the expression $A(A -\lambda \id)^{-1}(\phi)$ and proceeding as in the last expression, we can prove that 
%\bd
%(A-\lambda \id)^{-1} \colon \mathcal{H} \to D(A)
%\ed
%is continuous (where $D(A)$ is endowed with the graph norm).

In the same way we can prove that 
%$(A^*-\lambda \id)^{-1}$ is a continuous and negative linear operator and that
\bd
(A^*-\lambda \id)^{-1} \colon \mathcal{H} \to D(A^*)
\ed
is a continuous and negative linear operator and that and its explicit expression is given by
\bd
(A^*-\lambda \id)^{-1}(\phi)[r] = \frac{1}{\beta} \left ( -e^{\frac{\lambda}{\beta}r} \int_r^{\bar s} e^{-\frac{\lambda}{\beta}\tau} \phi(\tau) \ud \tau  \right )
\ed
Eventually we can define $B\nd 
%--- ((A-\lambda \id)(A^*-\lambda \id))^{-1}=  
%capire se questo pezzo serviva
(A^*-\lambda\id)^{-1} (A-\lambda\id)^{-1} = ((A-\lambda\id)^{-1})^* (A-\lambda\id)^{-1}$ that is continuous, positive and selfadjoint\footnote{See \cite{Yosida95} Proposition 2 page 273 for a proof of the equality $(A^*-\lambda\id)^{-1} = ((A-\lambda\id)^{-1})^*$.}. 
Moreover
\bd
(A^*-\lambda\id)B=(A-\lambda \id)^{-1} \leq 0
\ed
and so
\bd
A^*B= (A-\lambda \id)^{-1} + \lambda B  \leq \lambda B
\ed
if we choose $\lambda < 1$ we have that $A^*B$ is bounded and
\be
\label{eqRenardy}
A^*B \leq B
\ee
Thus $B$ satisfies all requirements of the \textquotedblleft weak case\textquotedblright \phantom{ } of \cite{CL4}.

\begin{Remark}
\label{Remark1}
Note that $B^{1/2}$ is a particular case of the operator that Renardy found in more generality in \cite{Renardy} and so $B^{1/2} \colon \mathcal{H} \to D(A^*)$ continuously and in particular $\mathcal{R}(B^{1/2}) \subseteq D(A^*)$. 
\end{Remark}

\begin{Notation}
For every $x\in\mathcal{H}$ we will indicate with $|x|_B$ the $B$-norm that is $\sqrt{\lla B  x,x \rra_{\mathcal{H}}}$. We will write $\mathcal{H}_{B}$ for the completion of $\mathcal{H}$ on the $B$-norm.
\end{Notation}

\begin{Remark}
Thanks to the definition of $A^*$ the graph norm on $D(A^*)$ is equivalent to the $H^1(0,\bar s)$ norm. In particular $D(A^*)$ is the the completion of 
\bd
K=\{ f|_{[0,\bar s]} \s : \s f\in C^{\infty}_c(\mathbb{R}) \s with \s supp(f) \subseteq (-\infty,\bar s) \}
\ed
with respect the $H^1(0,\bar s)$ norm.
So, since $H^1(0,\bar s) \hookrightarrow C([0,\bar s];\mathbb{R})$, we can calculate $\beta\delta_0$ on the elements of $D(A^*)$.
\end{Remark}

\begin{Notation}
\label{notationscalare}
In some cases the notation $\lla x,y\rra$ may be not clear, so when necessary we will use an index to avoid confusion: if $H$ is a Hilbert space (for example $H=\mathcal{H}\equiv L^2(0,\bar s)$ or $H=H^1(0,\bar s)$ or $D(A^*)$ ...) the notation $\lla x,y \rra_H$ will indicate the inner product in the Hilbert space $H$. Otherwise if $Z$ is an Banach space (possibly an Hilbert space) and $Z'$ its dual the notation $\lla x,y \rra_{Z\times Z'}$ will indicate the duality. In a few words, a single index means inner product, a double one indicates duality. 

When there is no index it is because we have omitted the index $\mathcal{H} \equiv L^2(0,\bar s)$.
\end{Notation}

\subsection{The control problem and the HJB equation}
In this subsection we describe the optimal control problem, state the hypotheses, define the HJB equation of the system and give a suitable definition of solution of the HJB equation.

We will consider the optimal control problem governed by the state equation
\be
\left\{ \begin{array}{ll}
\frac{\ud}{\ud s} i^* x(s) = A x(s) -\mu x(s) +\alpha(s) + \beta \delta_0 a(s)\\
x(0)=x
\end{array} \right.
\ee
that has a unique solution in the sense described in Section \ref{sezioneequazione}. Given two compact subsets $\Gamma$ and $\Lambda$ of $\mathbb{R}$ we consider the set of admissible boundary controls given by
\bd
\mathcal{A} \nd \left\{a\colon [0, +\infty) \to \Gamma\subseteq\mathbb{R} \s : \s a(\cdot) \s is \s measurable \right \}.
\ed
Moreover we call
\bd
\Sigma\nd \left\{ \gamma\colon [0,\bar s]\to \Lambda\subseteq\mathbb{R} \s : \s \gamma(\cdot)\in L^2(0,\bar s)   \right\}.
\ed
In view of the compactness of $\Lambda$ we have that $\Sigma\subseteq L^2(0,\bar s)$. We define the set of admissible distributed controls as
\bd
\mathcal{E} \nd \left \{\alpha\colon [0, +\infty) \to \Sigma\subseteq L^2(0,\bar s) \s : \s \alpha(\cdot) \s is \s measurable \right \}
\ed
In view of the complactness of $\Gamma$ and $\Lambda$ $\mathcal{A}\subseteq L^{2}_{loc}([0, +\infty);\mathbb{R})$ and $\mathcal{E}\subseteq L^{2}_{loc}([0, +\infty)\times[0,\bar s];\mathbb{R})$. We call $\|\Gamma \| \nd \sup_{a\in\Gamma} (|a|)$, $\|\Lambda \| \nd \sup_{b\in\Lambda} (|b|)$ and $\|\Sigma \| \nd \sup_{\alpha\in\Sigma} (|\alpha|_{\mathcal{H}=L^2(0,\bar s)})$ (they are bounded thanks to the boundedness of $\Gamma$ and $\Lambda$).

We will call \textit{admissible control} a couple $(\alpha(\cdot), a(\cdot)) \in \mathcal{E}\times \mathcal{A}$.
The cost functional will be of the form
\bd
J(x,\alpha(\cdot),a(\cdot))=\int_0^\infty e^{-\rho s} L(x(s),\alpha(s),a(s)) \ud s
\ed
where $L$ is uniformly continuous and satisfies the following conditions: there exists a $C_L\geq 0$ with
\bd
\begin{array}{ll}
(L1) \s\s |L(x,\alpha,a) - L(y,\alpha,a)| \leq C_L \lla B (x-y), (x-y) \rra_{\mathcal{H}\times\mathcal{H}} \forall (\alpha,a) \in \Sigma\times\Gamma\\
(L2) \s\s |L|\leq C_L <+\infty
\end{array}
\ed
We define formally the HJB equation of the system as
\be
\tag{HJB}
\label{HJB}
\rho u(x) - \lla \nabla u(x) , Ax \rra - \lla \nabla u (x), -\mu x \rra - H(x, \nabla u(x)) = 0
\ee
where $H$ is the Hamiltonian of the system and is defined as:
\bd
\left\lbrace 
\begin{array}{ll}
H\colon \mathcal{H} \times D(A^*) \to \mathbb{R}\\
H(x,p) \nd \inf_{(\alpha,a)\in\Sigma\times\Gamma} \left( \lla \beta\delta_0(p),a \rra_{\mathbb{R}} +\lla p,\alpha \rra_{\mathcal{H}} +  L(x,\alpha,a) \right) 
\end{array}
\right. 
\ed
(according to Notation \ref{notationscalare} $\lla \cdot , \cdot \rra_{\mathbb{R}}$ is the usual real product). 

Before we can introduce a suitable definition of (viscosity) solution of the HJB equation we have to give some preliminary definitions.
\begin{Definition}
\label{defblipschitz}
$v\in C(\mathcal{H})$ is \textit{Lipschitz with respect the $B$-norm} or  \textit{$B$-Lipschitz} if
there exists a constant $C$ such that $|v(x) - v(y)| \leq C |(x-y)|_B \nd C |B^{1/2}(x-y)|_{\mathcal{H}}$ for every choice of $x$ and $y$ in $\mathcal{H}$. In the same way we can give the definition of \textit{locally $B$-Lipschitz} function.
\end{Definition}
\begin{Definition}
A function $v\in C(\mathcal{H})$ is said to be $B$-continuous at a point $x\in \mathcal{H}$ if for every $x_n \in \mathcal{H}$ with $x_n \rightharpoonup x$ and $|B(x_n - x )| \rightarrow 0$, it holds that $v(x_n) \rightarrow v(x)$. In the same way we can define the $B$-upper/lower semicontinuity.
\end{Definition}
\begin{Definition}
We say that a function $\phi$ such that $\phi \in C^1(\mathcal{H})$ and $\phi$ is $B$-lower semicontinuous is a \textit{test function of type 1} and we will write $\phi \in test1$ if $\nabla \phi (x) \in D(A^*)$ for all $x\in\mathcal{H}$ and $A^* \nabla \phi\colon \mathcal{H} \to \mathcal{H}$ is continuous.
\end{Definition}
\begin{Definition}
We say that $g \in C^2(\mathcal{H})$ is a \textit{test function of type 2} and we will write $g \in test2$ if $g(x)=g_0(|x|)$ for some function $g_0\colon \mathbb{R}^+ \to \mathbb{R}$ nondecreasing.
\end{Definition}
\begin{Definition}
\label{defsubsol}
$u\in C(\mathcal{H})$ bounded and Lipschitz with respect the $B$-norm is a \textit{subsolution} of the HJB equation (or simply a \textquotedblleft \textit{subsolution}\textquotedblright) if for every $\phi\in test1$ and $g\in test2$ and a local maximum point $x$ of $u-(\phi+g)$ we have
\begin{multline}
\rho u(x) - \lla A^* \nabla \phi(x),x\rra - \lla \nabla \phi (x) + \nabla g(x) ,-\mu x \rra -\\
- \inf_{(\alpha, a ) \in \Sigma\times\Gamma}\Big ( \lla \beta\delta_0(\nabla\phi(x), a \rra_{\mathbb{R}} + \lla \nabla\phi(x)+\nabla g(x), \alpha \rra_{\mathcal{H}}+ L(x,\alpha,a) \Big ) \leq\\
\leq \frac{g_0'(|x|)}{|x|} \beta \frac{\|\Gamma\|^2}{2}
\end{multline}
\end{Definition}
\begin{Definition}
\label{defsupersol}
$v\in C(\mathcal{H})$ bounded and Lipschitz with respect the $B$-norm is a \textit{supersolution} of the HJB equation (or simply a "\textit{supersolution}") if for every $\phi\in test1$ and $g\in test2$ and a local minimum point $x$ of $v+(\phi+g)$ we have
\begin{multline}
\rho v(x) + \lla A^* \nabla \phi(x),x\rra + \lla \nabla \phi (x) + \nabla g(x) ,-\mu x \rra - \\
-\inf_{(\alpha, a ) \in \Sigma\times\Gamma}\Big ( - \lla \beta\delta_0(\nabla\phi(x), a \rra_{\mathbb{R}} - \lla \nabla\phi(x)+\nabla g(x), \alpha \rra_{\mathcal{H}} + L(x,\alpha,a) \Big )  \geq  \\
\geq - \frac{g_0'(|x|)}{|x|}\beta \frac{\|\Gamma\|^2}{2}
\end{multline}
\end{Definition}
\begin{Definition}
$v\in C(\mathcal{H})$ bounded and Lipschitz with respect the $B$-norm is a \textit{solution} of the HJB equation if it is at the same time a supersolution and a subsolution.
\end{Definition}
We write now some remarks on the definition we have just given, to explain its meaning:
\begin{Remark}
\label{remarkdefinizione}
In the definition of viscosity solution we have used two kinds of test function: the \textit{test1} and \textit{test2} that, as usual in the literature, play a different role in the definition. In view of their properties and their regularity the functions of the first set (test1) represent the \textquotedblleft good part\textquotedblright \phantom{ } while the main difficulties come from the function of the set test2 that have the role of localize the problem.

A difficulty of our case is the following: the trajectory is not Lipschitz in the norm of the Hilbert space $\mathcal{H}$ and so, given a function $g\in test2$, the term
\be
\label{eq10bis}
\frac{g(x(s)) - g(x)}{s}
\ee
(where $x(s)$ is a trajectory starting from $x$) cannot be treated with standard arguments. The idea is to consider only $B$-Lipschitz solution so that the maxima considered in Definition \ref{defsubsol} and Definition \ref{defsupersol} are in $D(A^*)$. If the \textquotedblleft starting point\textquotedblright \phantom{ } $x$ is in the domain of $D(A^*)$ there are some advantages in the estimate of (\ref{eq10bis}) but some problems remain: in such case we will prove in Proposition \ref{propg} that (if $\alpha(\cdot)$ is continuous)
\bd
\left | \frac{g(x(s)) - g(x)}{s} - \lla \nabla g (x), -\mu x + \alpha(0) \rra  \right | \leq \frac{g_0'(|x|)}{|x|} \beta \frac{\|\Gamma\|^2}{2} + O(s)
\ed
where the rest $O$ is uniform in the control. So the \textquotedblleft worse case\textquotedblright \phantom{ } is the one described in the definition.
\end{Remark}

\section{The value function and its properties}
\label{sectionvalore}
The value function is, as usual, the candidate to be the unique solution of the HJB equation. In this section we define the value function of the problem and then we verify that it has the regularity properties required to be a solution. Namely we will check that it is $B$-Lipschitz (Proposition \ref{propVBlipschitz}). To obtain such result we need to prove an approximation result (Proposition \ref{propU}) and then a suitable estimate for the solution of the state equation (Proposition \ref{propB-lip}).

The value function of our problem is defined as:
\bd
V(x) \nd \inf_{(\alpha(\cdot),a(\cdot)) \in \mathcal{E}\times \mathcal{A}} J(x, \alpha(\cdot), a(\cdot))
\ed
We consider the functions
\bd
\left \{
\begin{array}{l}
\eta_n\colon [0,\bar s] \to \mathbb{R}\\
\eta_n(r) \nd [2n - 2n^2 r]^+
\end{array}
\right .
\ed
(where $[\cdot]^+$ is the positive part). We then define
\bd
\left \{
\begin{array}{l}
\mathcal{C}_n^* \colon \mathbb{R} \to \mathcal{H}\\
\mathcal{C}_n^* \colon \gamma \mapsto \gamma \eta_n
\end{array}
\right .
\ed
Such functions are linear and continuous and their adjoints are
\be
\label{formadiCn}
\left \{
\begin{array}{l}
\mathcal{C}_n \colon \mathcal{H} \to \mathbb{R}\\
\mathcal{C}_n \colon x \mapsto \left\langle  x,\eta_n\right\rangle 
\end{array}
\right .
\ee
$\mathcal{C}_n^*$ \textquotedblleft approximate the delta measure\textquotedblright. The approximating equations we consider are
\be
\label{eqapproximatingstate}
\left\{ \begin{array}{ll}
\frac{\ud}{\ud s} x_n(s) = A x_n(s) -\mu x_n(s) + \alpha(s) + \beta\mathcal{C}_n^* a(s)\\
x_n(0)=x
\end{array} \right.
\ee
In the proofs we will use the mild solutions of the original and the approximating state equations. The first was introduced in (\ref{eqstatomildform}), the second can be found in (\cite{Pazy83} page 105 equation (2.3), we include the term $e^{-\mu s}$ in the semigroup):
\begin{multline}
x(s)= e^{- \mu s} e^{sA} x + \int_0^s e^{-(s-\tau)\mu} e^{(s-\tau)A} \alpha(\tau) \ud \tau-\\
 - A \int_0^s e^{-(s-\tau)\mu} e^{(s-\tau)A} ( a(\tau) \nu) \ud \tau
\end{multline}
\begin{multline}
x_n(s)= e^{- \mu s} e^{sA} x + \int_0^s e^{-(s-\tau)\mu} e^{(s-\tau)A} \alpha(\tau) \ud \tau+\\
+ \int_0^s e^{-(s-\tau)\mu} e^{(s-\tau)A} \beta \mathcal{C}_n^* a(\tau) \ud \tau
\end{multline}
\begin{Proposition}
\label{propU}
For $T>0$ and $(\alpha(\cdot),a(\cdot))\in\mathcal{E}\times\mathcal{A}$ 
\bd
\lim_{n\to \infty} \sup_{s \in [0,T]} |x_n(s) - x(s)|_{\mathcal{H}} =0
\ed
\end{Proposition}
\begin{proof}
Using the mild expressions we find
\begin{multline}
\left |  x(s) - x_n(s) \right | = \left | - A \int_0^s e^{-(s-\tau)\mu} e^{(s-\tau)A} (a(\tau) \nu) \ud \tau - \right .\\ 
\left . - \int_0^s e^{-(s-\tau)\mu} e^{(s-\tau)A} \beta \mathcal{C}_n^*(a(\tau)) \ud \tau \right |
\end{multline}
To estimate such expression we will use the explicit expression of the two terms (as two-variable function). We simplify the notation (only in this proof!) taking an \textquotedblleft extension\textquotedblright\phantom{ }of $a(\cdot)$ to the whole $\mathbb{R}$ obtained by putting $a(\cdot)$ identically $0$ on $\mathbb{R}^-$. So
\bd
y(s,r) \nd \left( - A \int_0^s e^{-(s-\tau)\mu} e^{(s-\tau)A} ( a(\tau)\nu) \ud \tau\right) [r] = e^{-\frac{\mu}{\beta}r} a(s-r/\beta)\\
\ed
\begin{multline}
y_n(s,r) \nd \left( \int_0^s e^{-(s-\tau)\mu} e^{(s - \tau) A}\beta \mathcal{C}_n^* (a(\tau)) \ud \tau \right ) [r]
= 
%\int_0^r e^{-\frac{\mu}{\beta}(r-\theta)} [2n-2n^2 \theta]^+ a \left(\frac{\theta-r}{\beta}+s\right) \ud \theta= 
\\=\int_0^{r\wedge (1/n)} e^{-\frac{\mu}{\beta}(r-\theta)} [2n-2n^2 \theta]^+ a \left (\frac{\theta-r}{\beta}+s \right ) \ud \theta
\end{multline}
Now for all $s\in [0,T]$
\begin{multline}
|y(s,\cdot) - y_n(s, \cdot)|^2_{\mathcal{H}=L^2(0,\bar s)} \leq \\ \leq \bigg ( \int_{1/n}^{\bar s} \left | e^{-\frac{\mu}{\beta}r} a(s-r/\beta) - \int_0^{1/n} e^{-\frac{\mu}{\beta}(r-\theta)} [2n-2n^2 \theta]^+ a \left (\frac{\theta-r}{\beta}+s \right ) \ud \theta \right |^2 \ud r \bigg ) + \\
+ \bigg ( \int_0^{1/n} \left | e^{-\frac{\mu}{\beta}r} a(s-r/\beta) - \int_0^{r} e^{-\frac{\mu}{\beta}(r-\theta)} [2n-2n^2 \theta]^+ a \left (\frac{\theta-r}{\beta}+s \right ) \ud \theta \right |^2 \ud r \bigg ) \leq
\end{multline}
(for $\bar s\leq T$)
\begin{multline}
\leq \Bigg ( e^{|\mu|s} \int_{0}^{T} \left | e^{-\mu (\frac{r}{\beta}-s)} a(s-r/\beta) - \right . \\
\left . - \int_0^{1/n} e^{-\mu(\frac{r-\theta}{\beta}-s)} [2n-2n^2 \theta]^+ a \left ( s+\frac{\theta-r}{\beta} \right ) \ud \theta \right |^2 \ud r \Bigg ) 
+ \bigg ( \frac{1}{n} e^{|\mu|/\beta T} 2 \|\Gamma \| \bigg ).
\end{multline}
Such estimate does not depends on $s$, the integral term goes to zero because it is the convolution of a function in $L^2(0,T)$ with an approximate unit and the second goes to zero for $n\to \infty$.
\end{proof}
\begin{Proposition}
\label{propS}
Let $\phi\in C^1(\mathcal{H})$ be such that $\nabla\phi \colon \mathcal{H} \to D(A^*)$ ($D(A^*)$ is endowed, as usual, with the graph norm) is continuous. Then, for an admissible control $(\alpha(\cdot),a(\cdot))$, if we call $x(\cdot)$ the trajectory starting from $x$ and subject to the control $(\alpha(\cdot),a(\cdot))$, we have that, for every $s>0$,
\begin{multline}
\label{phidixdis}
\phi(x(s)) = \phi(x) + \int_0^s \left[ \lla A^* \nabla \phi(x(\tau)), x(\tau) \rra + \lla \beta \delta_0(\nabla\phi(x(\tau))), a(\tau) \rra_{\mathbb{R}} + \right .\\
\left . + \lla \nabla \phi(x(\tau)), \alpha (\tau) \rra + \lla \nabla \phi(x(\tau)), - \mu x(\tau) \rra \right]  \ud \tau
\end{multline}
\end{Proposition}
\begin{proof}
In the approximating state equation (\ref{eqapproximatingstate}) the unbounded term $\beta\delta_0$ does not appear ($\beta\mathcal{C}_n^*$ are continuous) and then (see \cite{LiYong} Proposition 5.5 page 67) for every $\phi(\cdot)\in C^1(\mathcal{H})$ such that $A^*\nabla\phi(\cdot)\in C(\mathcal{H})$ we have
\begin{multline}
\phi(x_n(s)) = \phi(x) + \int_0^s \left[  \lla A^* \nabla \phi (x_n(\tau)) , x_n(\tau) \rra + \lla \nabla \phi(x_n(\tau)) , \beta \mathcal{C}_n^*  a(\tau) \rra + \right .\\
\left . + \lla \nabla \phi(x_n(\tau)) , \alpha(\tau) \rra +
\lla \nabla \phi(x_n(\tau)) , -\mu x_n(\tau) \rra
 \right ]\ud \tau.
\end{multline}
In view of the continuity of the operator $\mathcal{C}_n^*$ we can pass to its adjoint (see (\ref{formadiCn}) for an explicit form of the operator $\mathcal{C}_n$) and we obtain:
\begin{multline}
\label{phidixndis}
\phi(x_n(s)) = \phi(x) + \int_0^s \left[  \lla A^* \nabla \phi (x_n(\tau)) , x_n(\tau) \rra + \lla \beta \mathcal{C}_n \nabla \phi(x_n(\tau)) ,  a(\tau) \rra + \right .\\
\left . + \lla \nabla \phi(x_n(\tau)) , \alpha(\tau) \rra +
\lla \nabla \phi(x_n(\tau)) , -\mu x_n(\tau) \rra
 \right ]\ud \tau.
\end{multline}
Now we prove that every integral term of the (\ref{phidixndis}) converges to the corresponding term of the (\ref{phidixdis}). This fact, toghether with the pointwise convergence of ($\phi(x_n(s)) \xrightarrow{n \to \infty} \phi(x(s))$ due to Proposition \ref{propU}) prove the claim.

First we note that, in view of Proposition \ref{propU} and of the continuity of $x$, $x_n(\tau)$ is bounded uniformly in $n$ and $\tau\in[0,s]$ and, in view of the continuity of $\nabla\phi$, $\nabla\phi(x_n(r))$ is bounded uniformly in $n$ and $\tau\in[0,s]$
%\footnote{Note that here we use the compactness of the set $\{x(r) \; : \; r\in[0,s] \}$}. 
So we can apply the Lebesgue theorem (the pointwise convergence is given by Proposition \ref{propU} and $|\alpha(\tau)|\leq \|\Sigma\|$) and we prove that
\begin{multline}
\int_0^s \left[ \lla \nabla \phi(x_n(\tau)) , \alpha(\tau) \rra +
\lla \nabla \phi(x_n(\tau)) , -\mu x_n(\tau) \rra
 \right ]\ud \tau \xrightarrow{n\to\infty}\\
\xrightarrow{n\to\infty} \int_0^s \left[ \lla \nabla \phi(x(\tau)) , \alpha(\tau) \rra +
\lla \nabla \phi(x(\tau)) , -\mu x(\tau) \rra
 \right ]\ud \tau
\end{multline}
Now we observe that, in view of the continuity of $A^*\nabla\phi$ and of the of Proposition \ref{propU}, the term $A^* \nabla \phi (x_n(\tau))$ is bounded uniformly in $n$ and $\tau\in[0,s]$ so the same is true for
\[
| A^* \nabla \phi (x_n(\tau)) - A^* \nabla \phi (x(\tau)) |.
\]
Therefore we can use the Lebesgue theorem (the pointwise convergence is given by Proposition \ref{propU}) to conclude that
\bd
\int_0^s \lla A^* \nabla \phi (x_n(\tau)) , x_n(\tau) \rra \ud \tau \to \int_0^s \lla A^* \nabla \phi (x(\tau)) , x(\tau) \rra \ud \tau
\ed
We have now to prove that 
\be
\label{rimanedaprovare}
\int_0^s \lla \beta \mathcal{C}_n \nabla \phi(x_n(\tau)) , a(\tau) \rra \ud \tau \to \int_0^s \lla \beta \delta_0(\nabla\phi(x(\tau))), a(\tau) \rra_{\mathbb{R}} \ud \tau
\ee
We first note that $\mathcal{C}_n \xrightarrow{n \to \infty} \delta_0$ in $H^{-1} (0, \bar s)$ and then in $D(A^*)'$. Indeed given $z \in H^1(0, \bar s)$  we have
\begin{multline}
\left | (\mathcal{C}_n - \delta_0) z \right | = \left | \int_0^{\bar s} z[\tau] \eta_n[\tau] \ud \tau - z[0] \right | =\\
= \left | \int_0^{1/n} \left ( z[0] + \int_0^\tau \partial_\omega z[r] \ud r \right ) \eta_n[\tau] \ud \tau - z[0] \right |=
\end{multline}
($\partial_\omega z$ is the weak derivative of z) integrating by part
\begin{multline}
= \left| \left( z[0] + \int_0^{1/n} \partial_\omega z[r] \ud r \right ) \left ( \int_0^{1/n} \eta_n[r] \ud r \right ) - \right.\\
\left.- \int_0^{1/n} \partial_\omega z[\tau] \int_0^{\tau} \eta_n[r] \ud r \ud \tau - z[0] \right| \leq
\end{multline}
writing $\eta_n$ in explicit form and making calculi (note that $\int_0^{1/n} \eta_n[r] \ud r =1$)
\bd
\leq \left | \int_0^{\bar s} \chi_{[0,1/n]}[\tau] |\partial_\omega z[\tau]| \ud \tau \right | \leq \frac{1}{\sqrt{n}} \|z\|_{H^1(0,\bar s)}
\ed
Summarizing: by Proposition \ref{propU} $x_n(\cdot) \xrightarrow{n\to\infty} x(\cdot)$ in $C([0,T]; \mathcal{H})$, then (by hypothesis on $\phi$) $\nabla \phi (x_n(\cdot)) \xrightarrow{n\to\infty} \nabla \phi (x(\cdot))$ in $C([0,T]; D(A^*))$ and then, by the last estimate $\beta\mathcal{C}_n (\nabla \phi (x_n(\cdot))) \xrightarrow{n\to\infty} \beta\delta_0(\nabla \phi (x(\cdot)))$ in $C([0,T];\mathbb{R})$. Then (\ref{rimanedaprovare}) follows by Cauchy-Schwartz inequality (it is the scalar product in $L^2(0,s)$).
\end{proof}

\begin{Proposition}
\label{propB-lip}
Given $T>0$ and a control $(\alpha(\cdot),a(\cdot)) \in \mathcal{E}\times\mathcal{A}$ there exists $c_T$ such that for every $x,y \in \mathcal{H}$
\bd
\sup_{s\in[0,T]} |x_x(s) - x_y(s) |^2_B \leq c_T |x-y|^2_B,
\ed
where $x_y(\cdot)$ is the solution of
\bd
\left\{ \begin{array}{ll}
\frac{\ud}{\ud s} i^* x(s) = A x(s) +\alpha(s) -\mu x(s)  + \beta \delta_0 a(s)\\
x(0)=y
\end{array} \right.
\ed
and $x_x(\cdot)$ the solution with initial data $x$
\end{Proposition}
\begin{proof}
We use Proposition \ref{propS} with $\phi(x)=\lla Bx,x\rra$. So $\nabla \phi(x) = 2Bx$. We observe that $x_x(\cdot) - x_y(\cdot)$ satisfies the equation
\bd
\left\{ \begin{array}{ll}
\frac{\ud}{\ud s} i^*( x_x(s) - x_y(s))  = A (x_x(s)- x_y(s)) - \mu(x_x(s)- x_y(s)) \\
(x_x-x_y)(0)=x-y
\end{array} \right.
\ed
(the one of Proposition \ref{propS} with control identically $0$) and then by (\ref{eqRenardy})
\begin{multline}
|x_x(s) - x_y(s)|^2_B = |x-y|^2_B + 2 \int_0^s \lla A^*B (x_x(r)- x_y(r)) , (x_x(r)- x_y(r)) \rra -\\
-\mu \lla B(x_x(s)- x_y(s)), x_x(s)- x_y(s) \rra \ud r \leq
\end{multline}
\bd
\leq |x-y|^2_B + 2(1+|\mu|) \int_0^s \lla B (x_x(r)- x_y(r)) , (x_x(r)- x_y(r)) \rra \ud r
\ed
now we can use the Gronwall's lemma and obtain the claim.
\end{proof}

\begin{Proposition}
\label{propVBlipschitz}
Let $L$ satisfy (L1) and (L2). Then the value function $V$ is Lipschitz with respect the $B$-norm
\end{Proposition}
\begin{proof}
Assume $V(y) > V(x)$. Then we take $(\alpha(\cdot),a(\cdot))\in \mathcal{E}\times\mathcal{A}$ an $\varepsilon$-optimal control for $x$. We have:
\bd
|V(y)-V(x)| -\varepsilon \leq \int_0^{\infty} e^{-\rho t} |L(x_y(s),\alpha(s),a(s)) - L(x_x(s),\alpha(s),a(s)) | \ud s=
\ed
If we look the explicit for of $x_x(\cdot)$ and $x_y(\cdot)$ as two-variables functions we see that they depend on the initial data only for $s \in [0,\frac{\bar s}{\beta}]$. After this period they depends only on the control. So for $s>\frac{\bar s}{\beta}$ $x_x(s) = x_y(s)$ and so the previous integral is equal to
\bd
= \int_0^{\bar s/\beta} e^{-\rho t} |L(x_y(s),\alpha(s),a(s)) - L(x_x(s),,\alpha(s),a(s)) | \ud s \leq \s\s
\ed
(by (L1) and Proposition \ref{propB-lip})
\bd
\leq \int_0^{\bar s} e^{-\rho t} C_L |x_y(s)- x_x(s)|_B \ud s \leq \bar s c_{\bar s} C_L |x-y|_B
\ed
Letting $\varepsilon \to 0$ we have claim.
\end{proof}

\section{Existence and uniqueness of solution}
In this section we will prove that the value function is a viscosity solution of the HJB equation (Theorem \ref{thexistence}) and that the HJB equation admits at most one solution (Theorem \ref{thuniqueness}). \\

We remind that we use $\mathcal{H}_{B}$ to denote the completion of $\mathcal{H}$ in the $B$-norm. This notation will be used in the next propositions.

\begin{Proposition}
\label{prop0}
Let $u\in C(\mathcal{H})$ be a locally $B$-Lipschitz function. Let $\psi\in C^1(\mathcal{H})$, and let $x$ be a local maximum (or a local minimum) of $u-\psi$. Then $\nabla\psi(x) \in \mathcal{R}(B^{1/2}) \subseteq D(A^*)$.
\end{Proposition}
\begin{proof}
We do the proof only in the case in which $x$ is a local maximum (the other case is similar).

We take $\omega \in \mathcal{H}$ with $|\omega|=1$ and $h\in(0,1)$. Then for every $h$ small enough
\bd
\frac{(u(x-h\omega)-\psi(x-h\omega))}{h} \leq \frac{u(x) - \psi(x)}{h}
\ed
so 
\bd
\frac{\psi(x)-\psi(x-h\omega)}{h} \leq C |w|_B
\ed
and passing to the limit we have $\lla \nabla \psi(x), \omega \rra \leq C |\omega|_B$. Likewise
\bd
\frac{(u(x+h\omega)-\psi(x+h\omega))}{h} \leq \frac{u(x) - \psi(x)}{h}
\ed
so 
\bd
\frac{\psi(x)-\psi(x+h\omega)}{h} \leq C |w|_B
\ed
and passing to the limit we have $-\lla \nabla \psi(x), \omega \rra \leq C |\omega|_B$.

Putting together these two remarks we have
\bd
|\lla \nabla \psi (x) , \omega \rra| \leq C |\omega|_B
\ed
for all $\omega \in \mathcal{H}$.
So we can consider the linear extension of the continuous linear functional $\omega \mapsto \lla \nabla \psi (x) , \omega \rra$ to $\mathcal{H}_{B}$; we will call such extension $\Phi_x$ and by Riesz representation theorem we can find $z_x\in \mathcal{H}_B$ such that
\bd
\Phi_x(\omega) = \lla z_x, \omega \rra_{\mathcal{H}_{B}} \s\s \forall \omega\in\mathcal{H}_{B}
\ed
however
\begin{multline}
\lla z_x, \omega \rra_{\mathcal{H}_{B}} = \lla B^{1/2}(z_x), B^{1/2} (\omega) \rra_{\mathcal{H}} = \\ 
= \lla B^{1/2}(B^{1/2}(z_x)), \omega \rra_{(\mathcal{H}_{B})' \times (\mathcal{H}_{B})}
= \lla B^{1/2}(m_x), \omega \rra_{(\mathcal{H}_{B})' \times (\mathcal{H}_{B})}
\end{multline}
where $m_x \nd (B^{1/2}(z_x)) \in \mathcal{H}$. Now for $\omega \in \mathcal{H}$
\bd
\lla B^{1/2}(m_x), \omega \rra_{(\mathcal{H}_{B})' \times (\mathcal{H}_{B})} = \lla B^{1/2}(m_x), \omega \rra_\mathcal{H}
\ed
Therefore $\nabla \psi (x) = B^{1/2} (m_x) \in \mathcal{R}(B^{1/2})\subseteq D(A^*)$ where the last inclusion follows from Remark \ref{Remark1}.
\end{proof}
\subsection{Existence}
In this subsection we will prove that the value function is a solution of the HJB equation. In the next subsection we will prove that such solution is unique. We start with a lemma and two propositions. We will use the notation introduced in Remark \ref{notationquadre} on \textquotedblleft$x(s)$\textquotedblright\phantom{ } and \textquotedblleft$x[r]$\textquotedblright.
Moreover we will continue to use the symbol $\delta_0$ in the text so that $x[0]=\delta_0 x$ if $x\in D(A^*)$.

We have not found a simple reference for the following lemma so we prove it:
\begin{Lemma}
\label{lemmasobolev}
Let $x$ be a function of $H^1(0,\bar s)$ then
\begin{eqnarray}
& (i) & \lim_{s\to 0^+} \left(  \int_s^{\bar s} \frac{(x[r] - x[r-s])^2}{s} \ud r \right) =0\\
& (ii) & \lim_{s\to 0^+} \left(  \int_s^{\bar s -s} \frac{(x[r+s] - x[r])}{s} x[r] \ud r \right)= \frac{x^2[\bar s] - x^2[0]}{2}
\end{eqnarray}
\end{Lemma}
\begin{proof}
\textbf{part (i)}
\bd
\int_s^{\bar s} \frac{(x[r] - x[r-s])^2}{s} \ud r = \int_0^{\bar s} \psi_s[r] \ud r
\ed
where $\psi_s\colon [0, \bar s] \to \mathbb{R}$ is defined in the following way:
\bd
\psi_s[r]= \left\lbrace \begin{array}{ll}
			0 & if \s r\in [0,s)\\
			\frac{(x[r] - x[r-s])^2}{s} & if \s r\in[s, \bar s]
                        \end{array} \right. 
\ed
In order to prove the claim we want to apply the Lebesgue theorem. First we will see the a.e. convergence of the $\psi_s$ to zero: for $r>0$ we take $s<r$:
\bd
\psi_s[r]\leq \frac{\left | \int_{r-s}^{r} \partial_\omega x[\tau] \ud \tau \right |}{s} \left | x[r] - x[r-s] \right|
\ed
where $\partial_\omega x$ is the weak derivative of $x$ ($x$ is in $H^1$ for hypothesis). Now almost every $r$ is a Lebesgue point and then
\bd
\frac{\left | \int_{r-s}^{r} \partial_\omega x(\tau) \ud \tau \right |}{s} \xrightarrow{s\to 0^+} |\partial_\omega x[r]|\s\s\s\s\s a.e.\s in\s r\in(0,\bar s]
\ed
while the part $\left | x[r] - x[r-s] \right|$ goes uniformly to $0$.

In order to dominate the convergence we note that by Morrey's theorem (\cite{Evans98} Theorem 4 page 266) every $x\in H^1(0,\bar s)$ is $1/2$-Holder then there exists a positive $C$ such that for every $s\in (0,\bar s]$ and every $r\in [s, \bar s]$ we have
\bd
\frac{|x[r] - x[r-s]|}{\sqrt{s}}\leq C
\ed
and then
\bd
\frac{|x[r] - x[r-s]|^2}{s}\leq C^2
\ed
this allows to dominate $\psi_s$ with the constant $C^2$, use the Lebesgue theorem and obtain the claim.\\
\textbf{part (ii):} 
\begin{multline}
I(s) \nd \int_s^{\bar s -s} \frac{(x[r+s] - x[r])}{s} x[r] \ud r = \\
= \int_s^{\bar s -s} \frac{(x[r+s] x[r])}{s} \ud r - \int_0^{\bar s - 2s} \frac{(x[r+s] x[r+s])}{s} \ud r = \\
= -\int_s^{\bar s -2s} \frac{(x[r+s] - x[r])}{s} x[r+s] \ud r + \int_{\bar s - 2s}^{\bar s -s} \frac{(x[r+s] x[r])}{s} \ud r + \\
+ \int_0^{s} - \frac{(x[r+s])^2}{s} \ud r \nd - I_1(s) + I_2(s) + I_3(s)
\end{multline}
By the continuity of $x$ we see that:
\bd
I_2(s) \xrightarrow{s\to 0^+} x^2[\bar s]
\ed
and
\bd
I_3(s) \xrightarrow{s\to 0^+} - x^2[0]
\ed
Moreover, using similar arguments that in (i) we find that
\begin{multline}
\lim_{s\to 0^+} (I (s) - I_1(s)) = \lim_{s\to 0^+} \int_s^{\bar s - 2s} - \frac{(x[r+s] - x[r])^2}{s} \ud r +\\
+ \lim_{s\to 0^+} \int_{\bar s - 2s}^{\bar s -s} \frac{(x[r+s]- x[r])}{s} x[r] \ud r = 0
\end{multline}
so the limit $\lim_{s\to 0^+} I(s)$ exist if and only if there exist the limit $\lim_{s \to 0^+} \frac{I_1(s) + I(s)}{2}$ and in such case they have the same value. But
\bd
\frac{I_1(s) + I(s)}{2} = \frac{I_2(s) + I_3(s)}{2} \xrightarrow{s\to 0^+} \frac{x^2[\bar s] - x^2[0]}{2}
\ed
and then $\lim_{s\to 0^+} \left(  \int_s^{\bar s -s} \frac{(x[r+s] - x[r])}{s} x[r] \ud r \right)= \frac{x^2[\bar s] - x^2[0]}{2}$.
\end{proof}

\begin{Lemma}
\label{lemmauniforme}
Given $x\in D(A^*)$ there exists a real function $O(s)$ such that $O(s)\xrightarrow{s\to 0} 0$ and such that for every control $(\alpha(\cdot),a(\cdot))\in\mathcal{E}\times\mathcal{A}$ we have that
\[
\left | x(s) -x \right | \leq O(s)
\]
(where we called $x(s)$ the trajectory that starts from $x$  and subject to the control $(\alpha(\cdot),a(\cdot))$). Note that $O(s)$ is independent of the control.
\end{Lemma}
\begin{proof}
We consider $s\in(0,1]$. This is an arbitrary choice but we are interested only in the behavior of $x(\cdot)$ near to $0$ so we can assume it without problems.
We use the explicit expression of $x(s,r)$:
\begin{multline}
\| x(s) - x \|^2_{\mathcal{H}=L^2(0,\bar s)} = \\
=\int_{\beta s}^{\bar s} \left | e^{-\mu s} x[r-\beta s] + \int_0^s e^{-\mu \tau} \alpha(s-\tau, r-\beta \tau) \ud \tau - x[r] \right |^2 \ud r + \\
+ \int_0^{\beta s} \left | e^{-\frac{\mu}{\beta}r} a(s-r/\beta) + \int_0^{r/\beta} e^{-\mu \tau} a(s-\tau, r- \beta\tau) \ud \tau -x[r] \right |^2 \ud r \leq\\
\leq 2 \int_{\beta s}^{\bar s} \left | e^{-\mu s} x[r-\beta s] - x[r] \right |^2 \ud r + 2 \int_{\beta s}^{\bar s} \left | \int_0^s e^{|\mu|} \|\Gamma \| \ud \tau \right |^2 \ud r +\\
+ \int_0^{\beta s} \left | e^{|\mu|} \|\Gamma \| + \int_0^{r/\beta} e^{|\mu|} \| \Lambda \| \ud \tau + |x|_{L^\infty(0,\bar s)} \right |^2 \ud r \leq
\end{multline}
(We have used that $x\in D(A^*)\subseteq W^{1,2}(0,\bar s)$ so it is continuous and $|x|_{L^\infty(0,\bar s)}<+\infty$)
\begin{multline}
\leq 2 \int_0^{\bar s} \left | e^{-\mu s} x[(r-\beta s)\wedge 0] - x[r] \right |^2 \ud r + 2 s^2 \bar s \left (e^{|\mu|} \|\Gamma\| \right )^2 +\\ 
+ s \beta \left (e^{|\mu|} \|\Gamma\| + |x|_{L^\infty} + s e^{|\mu|} \|\Lambda\| \right )^2 
\end{multline}
Observe that in this estimate the control $(\alpha(\cdot),a(\cdot))$ does not appear. The second and the third terms goes to zero for $s\to 0$. In the first we can use Lebesgue theorem observing that
\[
\left | e^{-\mu s} x[(r-\beta s)\wedge 0] - x[r] \right | \leq \left ( e^{|\mu|} |x|_{L^\infty} + |x|_{L^\infty}  \right ) \;\;\; \forall (s,r)\in (0,1]\times [0,\bar s]
\]
and that $\left | e^{-\mu s} x[(r-\beta s)\wedge 0] - x[r] \right | \xrightarrow{s\to 0}0$ pointwise. So the statement is proven.
\end{proof}

\begin{Proposition}
\label{propg}
Given $x\in D(A^*)$ and $g\in test2$ there exists a real function $O(s)$ such that $O(s)\xrightarrow{s\to 0} 0$ and such that for every control $(\alpha(\cdot),a(\cdot))\in\mathcal{E}\times\mathcal{A}$  with $a(\cdot)$ continuous we have that
\bd
\left | \frac{g(x(s)) - g(x)}{s} - \frac{\int_0^s \lla \nabla g(x) , \alpha(r) \rra}{s} - \lla \nabla g(x) , -\mu x \rra \right | \leq \frac{g_0'(|x|)}{|x|} \beta \frac{\|\Gamma\|^2}{2} + O(s)
\ed
(where we called $x(s)$ the trajectory that starts from $x$  and subject to the control $(\alpha(\cdot),a(\cdot))$). Note that $O(s)$ is independent of the control.
\end{Proposition}
\begin{proof}
First we write
\begin{multline}
\label{eqgxs-gx}
\frac{g(x(s)) - g(x)}{s} - \lla \nabla g(x) , -\mu x\rra - \frac{\int_0^s \lla \nabla g(x), \alpha(r) \rra}{s} = \\
= \frac{g(x(s)) - g(y(s)) + g(y(s)) - g(x)}{s} - \lla \nabla g(x) , -\mu x\rra - \frac{\int_0^s \lla \nabla g(x), \alpha(r) \rra}{s}
\end{multline}
where $y(\cdot)$ is the solution of
\begin{equation}
\label{eqdefdiy}
\left\lbrace 
\begin{array}{l}
\dot{y}(s)= A y(s) + \beta \delta_0 a(s)\\
y(0)=x
\end{array}
\right .
\end{equation}
(that is our system when $\mu=0$ and $\alpha(\cdot)=0$). $x(\cdot)$ satisfies the mild equation\footnote{We have already written an explicit mild form of the solution in (\ref{eqstatomildform}), the form we use here is different, indeed it is not explicit because the $x$ appears also in the second term. The only difference between the two formula is the following: equation (\ref{eqstatomildform}) is the equation we obtain if we include the term $-\mu x$ in the generator of the semigroup, equation (\ref{eqstatomildformdiversa}) is the form we obtain if we maintain the term $-\mu x$ out of the generator of the semigroup. The two forms are equivalent.}
\begin{equation}
\label{eqstatomildformdiversa}
x(s)=e^{sA}x - A \int_0^s e^{(s-\tau)A} (a(\tau)\nu) \ud \tau + \int_0^s e^{(s-\tau)A} (\alpha(\tau)-\mu x(\tau))\ud \tau
\end{equation}
The term $\left ( e^{sA}x - A \int_0^s e^{(s-\tau)A} (a(\tau)\nu) \ud \tau \right )$ is the mild solution of $y(\cdot)$ and
\[
x(s) - y(s) = \int_0^s e^{(s-\tau)A} (\alpha(\tau)-\mu x(\tau))\ud \tau.
\]
Now we come back to (\ref{eqgxs-gx}), we have
\begin{multline}
\label{eqfu43}
\left | \frac{g(x(s)) - g(x)}{s} - \frac{\int_0^s \lla \nabla g(x), \alpha(r) \rra \ud r}{s} -\lla \nabla g(x), -\mu x\rra \right |\leq \\
\leq
\left | \frac{g(x(s)) - g(y(s))}{s} - \frac{\int_0^s \lla \nabla g(x), \alpha(r) \rra \ud r}{s} -\lla \nabla g(x), -\mu x\rra \right |
+ \left | \frac{g(y(s)) - g(x)}{s} \right |.
\end{multline}
In order to estimate the first addendum we use the Taylor expansion as follows:
\begin{multline}
\frac{g(x(s)) - g(y(s))}{s} = \lla \nabla g(y(s)) , \frac{x(s) - y(s)}{s} \rra + \\
+ \lla \nabla g (\xi(s)) - \nabla g(y(s)) , \frac{x(s) - y(s)}{s} \rra=
\end{multline}
where $\xi(s)$ is a point between $x(s)$ and $y(s)$
\begin{multline}
\label{eqfu45}
= \lla \nabla g(y(s)) , \frac{\int_0^s e^{(s-\tau)A} (\alpha(\tau)-\mu x(\tau))\ud \tau}{s} \rra + \\
+ \lla \nabla g (\xi(s)) - \nabla g(y(s)) , \frac{\int_0^s e^{(s-\tau)A} (\alpha(\tau)-\mu x(\tau))\ud \tau}{s} \rra
\end{multline}
We know by Lemma \ref{lemmauniforme} that $x(s)\xrightarrow{s\to 0}x$ $y(s)\xrightarrow{s\to 0}x$ uniformly in the control $(\alpha(\cdot), a(\cdot))$, and so $\nabla g(y(s)) \xrightarrow{s\to 0}\nabla g (x)$ uniformly in the control and $\left | \nabla g(y(s)) - \nabla g (\xi(s)) \right | \xrightarrow{s\to 0}0$ uniformly in the control. Moreover, in view of boundedness of the control and of the fact that $x(s)\xrightarrow{s\to 0}x$ uniformly in the control (Lemma \ref{lemmauniforme}) we can prove that the term
\[
\left |\frac{\int_0^s e^{(s-\tau)A} (\alpha(\tau)-\mu x(\tau))\ud \tau}{s} \right|_{\mathcal{H}}
\]
is bounded uniformly in the control and $s\in(0,\bar s]$ and we conclude that the second term of the (\ref{eqfu45}) goes to zero uniformly in $(\alpha(\cdot),a(\cdot))$ and that
\begin{equation}
\label{eqstimagxs-gys}
\left | \frac{g(x(s)) - g(y(s))}{s} - \frac{\int_0^s \lla \nabla g(x), \alpha(r) \rra \ud r}{s} -\lla \nabla g(x), -\mu x\rra \right | \leq O(s)
\end{equation}
where $O(s)\xrightarrow{s\to 0}0$ and it does not depend on the control.

So we have now to estimate the second term of the (\ref{eqfu43}), namely $\left | \frac{g(y(s)) - g(x)}{s} \right |$. If we prove that it is smaller then $\frac{g_0'(|x|)}{|x|} \beta \frac{\|\Gamma\|^2}{2} + O(s)$ where $O(s)$ does not depend on the control we have proven the proposition.

We first note that
\bd
\nabla g(x) = g_0'(|x|) \frac{x}{|x|}\s
\ed
and
\bd 
D^2 g(x) = g_0''(|x|) \frac{x}{|x|}\otimes\frac{x}{|x|} + g_0'(|x|) \left (
\frac{\id}{|x|} - \frac{x\otimes x}{|x|^3} \right )
\ed
We consider the Taylor's expansion of $g$ at $x$:
\begin{multline}
\frac{g(y(s)) - g(x)}{s} = \frac{\lla \nabla g (x), y(s) - x \rra }{s} + \frac{1}{2} \frac{(y(s)-x)^T (D^2 g(x)) (y(s)-x)}{s} + \\
+ \frac{o (|y(s) -x|^2)}{s}=\\
=
\frac{g_0'(|x|)}{|x|} \left ( \lla x, \frac{y(s) - x}{s} \rra + \frac{1}{2} \frac{\lla y(s)-x, y(s) - x \rra}{s} \right ) + \\
+ \frac{1}{2} \left ( \frac{g_0''(|x|)}{|x|^2} - \frac{g_0'(|x|)}{|x|^3}   \right ) \frac{\lla x, y(s) - x \rra^2}{s} + \frac{o (|y(s) -x|^2)}{s} \nd \\
\nd \s P1 \s + \s P2 \s + \s P3 
\end{multline}
First we prove that $P2$ and $P3$ go to zero uniformly in $(\alpha(\cdot),a(\cdot))$ and then we will estimate $P1$. We proceed in two steps:\\
\textbf{step 1}: \textit{There exists a constant $C$ such that for every admissible control $(\alpha(\cdot),a(\cdot))\in\mathcal{E}\times\mathcal{A}$  with $a(\cdot)$ continuous and every $s\in (0,1]$\footnote{In the expression of $y(\cdot)$ the distributed control $\alpha(\cdot)$ does not appear, so we will speak from now only of the boundary control $a(\cdot)$}
\bd
\left |\frac{\lla x, y(s) - x \rra}{s} \right | \leq C
\ed
}
(as before the choice of the interval (0,1] it is not essential: we are interested in the behavior near zero).
We observe first that the explici solution of $y(s)[r]$ can be found taking $\mu=0$ and $\alpha=0$ in (\ref{eq5bis}). We have:
\[
y(s,r)=
\left\{ \begin{array}{ll}
x(r-\beta s) & \; r\in[\beta s, \bar s]\\
a(s-r/\beta) & \; r\in [0,\beta s)
\end{array} \right.
\]
so
\begin{multline}
\label{eqdistep1}
\frac{\lla x, y(s) - x \rra}{s} = \int_{\beta s}^{\bar s} x[r] \frac{(x[r-\beta s] - x[r])}{s} \ud r + \frac{\int_0^{\beta s} x[r] ( a(s-r/\beta) - x[r]) \ud r}{s} = \\
= \int_{\beta s}^{\bar s - \beta s} x[r] \frac{(x[r+\beta s] - x[r])}{s} \ud r + \frac{\int_{\bar s -\beta s}^{\bar s} - x^2[r] \ud r}{s} + \frac{\int_0^{\beta s} x[r] x[r+\beta s] \ud r}{s} + \\
+ \frac{\int_0^{\beta s} x[r] a(s-r/\beta) \ud r}{s} - \frac{\int_0^{\beta s} x^2[r] \ud r}{s}
\end{multline}
The third and the fifth part have opposite limits, the second goes to zero thanks to the fact that $x \in D(A^*)$ and then $x$ is continuous and $x(\bar s)=0$. The first part goes to $ -\frac{\beta}{2}x^2[0]=\lla A^* x, x \rra$ in view of Lemma \ref{lemmasobolev}. The only term in which the control appears is the fourth but we can estimate it as follows:
\bd
\left |\frac{\int_0^{\beta s} x[r] a(s-r/\beta) \ud r}{s} \right | \leq \frac{\int_0^{\beta s} |x[r]| \|\Gamma\| \ud r}{s} \leq \beta \max_{r\in [0, \bar s]} |x[r]| \|\Gamma\|
\ed

\textbf{step 2}: \textit {There exists a constant $C$ such that for every admissible control $a(\cdot)\in\mathcal{A}$  with $a(\cdot)$ continuous and every $s\in (0,1]$
\bd
\frac{|y(s)-x|^2}{s} \leq C
\ed
}

Indeed
\begin{multline}
\label{eqsecondascomposizione}
\left | \frac{\lla y(s)-x , y(s) - x \rra}{s} \right | = \left | \int_{\beta s}^{\bar s} \frac{(x[r-\beta s]-x[r])^2}{s} \ud r \right | +\\
+\left | \frac{\int_0^{\beta s} (a(s-r/\beta)-x[r])^2 \ud r}{s}\right |
\end{multline}
in view of the fact that $x\in D(A^*) \subseteq H^1(0,\bar s)$ and of the Lemma \ref{lemmasobolev} the first part goes to zero. Moreover, since $x\in H(0,\bar s)\subseteq L^{\infty}(0,\bar s)$, the second part is less or equal to
\be
\label{eqterzascomposizione}
\frac{\int_0^{\beta s} \|\Gamma\|^2 \ud r}{s} + \frac{\int_0^{\beta s} 2 |x[r]| \|\Gamma\| \ud r}{s} + \frac{\int_0^{\beta s} |x[r]|^2 \ud r}{s} \leq C.
\ee
This completes step 2.

From step 2 it follows that 
\bd
\frac{o(|y(s) -x|^2)}{s}=\frac{o(|y(s) -x|^2)}{|y(s) -x|^2} \frac{|y(s) -x|^2}{s} \xrightarrow{s\to 0^+} 0
\ed
uniformly in $a(\cdot)$. Thus $|P3| \xrightarrow{s\to 0} 0$ uniformly in $a(\cdot)$. Moreover
\bd
\frac{\lla x, y(s) - x \rra^2}{s} \leq \frac{|\lla x, y(s) - x \rra|}{s} |x| |y(s)-x|
\ed
and so, from step 1 and Lemma \ref{lemmauniforme}, $|P2| \xrightarrow{s\to 0} 0$ uniformly in $a(\cdot)$.\\

\textbf{step 3}: \textit {Conclusion}

We now estimate $P1$. We can write a more explicit form of $P1$ as in the proofs of step 1 and step 2 ((\ref{eqdistep1}), (\ref{eqsecondascomposizione}) and (\ref{eqterzascomposizione})) and using the same arguments we can see that there exists a rest $o(1)$ (depending only on $x$) with $o(1) \xrightarrow{s\to 0} 0$ such that for every control $a(\cdot)$ continuous
\begin{multline}
P1=\frac{g_0'(|x|)}{|x|} \left ( \lla A^* x,x \rra + \frac{\int_0^{\beta s} x[s] a(s-r/\beta) \ud r}{s} + \frac{1}{2} \frac{\int_0^{\beta s} (a(s-r/\beta))^2 \ud r}{s} + \right.\\
\left. + \frac{1}{2} \frac{\int_0^{\beta s} x^2[r] \ud r}{s} + \frac{1}{2} \frac{\int_0^{\beta s} - 2 x[r] a(s-r/\beta) \ud r}{s} \right ) + o(1)
\end{multline}
The fourth part of the above, that does not depend on the control, goes to $\beta\frac{x[0]^2}{2}$ that is the opposite of the first part. The second and the fifth part are opposite. So we have that
\bd
P1= o(1) + \frac{g_0'(|x|)}{|x|} \left (\frac{1}{2}\frac{\int_0^{\beta s} (a(s-r/\beta))^2 \ud r}{s} \right ) \leq o(1) + \frac{1}{2}\frac{g_0'(|x|)}{|x|} \beta \|\Gamma \|^2
\ed
Now, using the estimates on $P1$, $P2$ and $P3$ we see that
\[
\left | \frac{g(y(s))-g(x)}{s} \right | \leq O(s) + \frac{1}{2}\frac{g_0'(|x|)}{|x|}\beta \|\Gamma\|^2.
\]
Using this fact and equation (\ref{eqstimagxs-gys}) in (\ref{eqfu43}) we have proven the proposition.
\end{proof}

\begin{Proposition}
\label{propphi}
If $x\in D(A^*)$ and $\phi\in test1$ then there exists a real function $O(s)$ such that $O(s)\xrightarrow{s\to 0} 0$ and such that for every control $(\alpha(\cdot),a(\cdot))\in\mathcal{E}\times\mathcal{A}$  with $a(\cdot)$ continuous we have that
\begin{multline}
\left | \frac{\phi(x(s)) - \phi(x)}{s} - \frac{\int_0^s \lla \nabla \phi (x), \alpha(r) \rra \ud r}{s} - \lla \nabla \phi (x), - \mu x \rra - \right . \\ 
\left . - \lla A^* \nabla\phi(x), x \rra - \frac{\int_0^s \lla \beta \delta_0(\nabla\phi(x)), a(r) \rra_{\mathbb{R}} \ud r}{s} \right | \leq O(s)
\end{multline}
(where we called $x(s)$ the trajectory that starts from $x$  and subject to the control $(\alpha(\cdot),a(\cdot))$). Note that $O(s)$ is independent of the control.
\end{Proposition}
\begin{proof}
We proceed as in the proof of Proposition \ref{propg} observing that 
\[
\frac{\phi(x(s)) - \phi(x)}{s} = \frac{\phi (x(s)) - \phi (y(s))}{s} + \frac{\phi(y(s)) - \phi(x)}{s}
\]
where $y(\cdot)$ is the solution of (\ref{eqdefdiy}). It is possible to prove, using exactly the same arguments used in the proof of Proposition \ref{propg} that
\[
\left | \frac{\phi(x(s)) - \phi (y(s))}{s} - \lla \nabla \phi (x), -\mu x \rra - \frac{\int_0^s \lla \nabla\phi(x), \alpha(r)\rra \ud r }{s}  \right |\leq O(s)
\]
where $O(s)\xrightarrow{s\to 0}0$ and does not depend on the control. So we have to prove that
\[
\left | \frac{\phi(y(s)) - \phi (x)}{s} - \lla A^*\nabla \phi (x), x \rra - \frac{\int_0^s \beta \lla \delta_0 \nabla\phi(x), a(r)\rra_{\mathbb{R}} \ud r } {s}  \right |\leq O(s)
\]
where $O(s)\xrightarrow{s\to 0}0$ and does not depend on the control.

We write
\begin{multline}
\frac{\phi(y(s))-\phi(x)}{s} = I_0 + I_1 \nd \lla \nabla \phi(x), \frac{y(s)-x}{s} \rra + \\
+ \lla \nabla \phi (\xi(s)) - \nabla \phi (x) , \frac{y(s)-x}{s} \rra
\end{multline}
where $\xi(s)$ is a point between $x$ and $y(s)$. In view of Lemma \ref{lemmauniforme}, $|y(s)-x| \xrightarrow{s\to 0} 0$ uniformly in the control, so $|\xi(s)-x| \xrightarrow{s\to 0} 0$ uniformly in $a(\cdot)$. By hypothesis
\bd
\nabla \phi \colon \mathcal{H} \to D(A^*) \s\s\s and\s it \s is \s continuous 
\ed
($D(A^*)$ is endowed with the graph norm). Then
\be
\label{eqphi1}
|\nabla \phi (\xi(s)) - \nabla \phi(x) |_{D(A^*)} \xrightarrow{s\to 0} 0
\ee
uniformly in $a(\cdot)$.\\

If we read equation (\ref{eqdefdiy}) in $D(A^*)'$ it appears as an equation of the form
\bd
\left\lbrace 
\begin{array}{ll}
\dot u (t) = \tilde A u (t) + f(t)\\
u(0)=x
\end{array}
\right.
\ed
where $f(t)$ is a bounded measurable function ($|f(t)|_{D(A^*)'} \leq \beta |\delta_0|_{D(A^*)'} \|\Gamma\|$) and $\tilde A$ is an extension of $A$ that generates of a $C_0$-semigroup on $D(A^*)'$. So\footnote{In view of the fact that $x$ is in $\mathcal{H}\subseteq D(\tilde A)\subseteq D(A^*)'$, see \cite{Silviatesi} for a proof.} we can choose a constant $C$ that depends on $x$ such that, for all admissible control $a(\cdot)$ continuous and all $s\in (0,1]$,
\be
\label{eqphi2}
\frac{|y(s) - x|_{D(A^*)'}}{s} \leq C
\ee
Thus by (\ref{eqphi1}) and (\ref{eqphi2}), we can say that $|I_1| \xrightarrow{s\to 0} 0$ uniformly in $a(\cdot)$. Therefore
\bd
\left | \frac{\phi(y(s))-\phi(x)}{s} - \frac{\lla \nabla \phi(x), y(s) -x \rra}{s} \right | \xrightarrow{s\to 0} 0
\ed 
uniformly in $a(\cdot)$. We now write
\begin{multline}
\frac{\lla \nabla \phi(x), y(s) -x \rra}{s} = \int_{\beta s}^{\bar s} \nabla \phi (x)[r] \frac{(x[r-\beta s]-x[r])}{s} \ud r +\\
+ \frac{\int_0^{\beta s} \nabla \phi(x)[r] (a(s-r/\beta) - x[r]) \ud r}{s} = \\
= \int_{\beta s}^{\bar s -\beta s} x[r] \frac{\nabla\phi(x)[r+\beta s]-\nabla \phi (x)[r]}{s} \ud r + \int_{\bar s -\beta s}^{\bar s} \frac{(-\nabla\phi(x)[r]x[r])}{s} \ud r +\\ 
+ \frac{ \int_{0}^{\beta s} (\nabla\phi(x)[r+\beta s]x[r]) \ud r }{s}  + \frac{\int_0^{\beta s} \nabla \phi(x)[r] a(s-r/\beta) \ud r}{s} +\\ 
+ \frac{\int_0^{\beta s} -\nabla \phi(x)[r] x[r] \ud r}{s} 
\end{multline}
The third and the fifth terms, that do not depend on the control, have opposite limits, the second goes to zero because $\nabla\phi(x)$ and $x$ are in $D(A^*)$ and then $x[\bar s] = 0 = \nabla \phi (x) [\bar s]$. The first term goes to $\lla A^* \nabla\phi(x),x\rra$. Finally we observe that the only term that depends on the control is the fourth and
\bd
\left | \frac{\int_0^{\beta s} \nabla \phi(x)[r] a(s-r/\beta) \ud r}{s} - \beta \frac{\int_0^s \nabla \phi(x)[0] a(s-r') \ud r'}{s} \right | \xrightarrow{s\to 0} 0
\ed
uniformly in $a(\cdot)$ and, since $\phi(x)[0]$ is a constant,
\bd
\beta \frac{\int_0^s \nabla \phi(x)[0] a(s-r) \ud r}{s}=\frac{\int_0^s \lla \beta \delta_0 \nabla \phi(x), a(r) \rra_{\mathbb{R}} \ud r}{s}
\ed
This complete the proof.
\end{proof}

We can now prove that the value function is a solution of the HJB equation equation.
\begin{Theorem}
\label{thexistence}
Let $L$ satisfy (L1) and (L2) let $\Gamma$ and $\Lambda$ be a compact subsets of $\mathbb{R}$. Then the value function $V$ is bounded, $B$-Lipschitz and is a solution of the HJB equation.
\end{Theorem}
\begin{proof}
The boundedness of $V$ follows from the boundedness of $L$ (assumption (L2)). The $B$-Lipschitz property is the result of Proposition \ref{propVBlipschitz}. It remains to verify that $V$ is a solution of the HJB equation.\\
\textbf{Subsolution:}\\
Let $x$ be a local maximum of $V-(\phi+g)$ for $\phi \in test1$ and $g\in test2$. Thanks to Proposition \ref{prop0} we know that $\nabla(\phi+g)(x)\in D(A^*)$. Moreover we know that $\nabla \phi (x) \in D(A^*)$ for the definition of the set test1. So $\nabla g(x) = g_0'(|x|) \frac{x}{|x|} \in D(A^*)$ and this implies that $x\in D(A^*)$. We can assume that $V(x)-(\phi+g)(x)=0$. We consider the constant control $(\alpha(\cdot),a(\cdot))\equiv (\alpha,a)\in\Sigma\times\Gamma$ and $x(s)$ the trajectory starting from $x$ and subject to $(\alpha,a)$.
Then for $s$ small enough
\bd
V(x(s))- (\phi+g)(x(s)) \leq V(x) - (\phi+g)(x)
\ed
and thanks to the Bellman principle of optimality we know that
\bd
V(x) \leq e^{-\rho s} V(x(s)) + \int_0^s e^{-\rho r} L(x(r),\alpha,a) \ud r
\ed
Then
\begin{multline}
\frac{1-e^{-\rho s}}{s} V(x(s)) - \frac{\phi(x(s)) - \phi(x)}{s} - \frac{g(x(s)) - g(x)}{s} - \\
- \frac{\int_0^s e^{-\rho r} L(x(r),\alpha,a) \ud r}{s} \leq 0.
\end{multline}
Using Proposition \ref{propg} and Proposition \ref{propphi} we can now pass to the limsup as $s\to 0$ to obtain
\begin{multline}
\rho V(x) - \lla \nabla\phi(x), -\mu x \rra - \lla \nabla g(x), -\mu x \rra - \\
- \bigg ( \lla A^* \nabla\phi(x),x\rra + \lla \beta\delta_0(\nabla\phi(x)),a\rra_{\mathbb{R}} + \lla \nabla\phi(x) , \alpha \rra + \lla \nabla g(x) , \alpha \rra + L(x,\alpha,a)  \bigg ) \leq\\
 \leq \frac{g_0'(|x|)}{|x|}\beta\frac{\|\Gamma\|^2}{2}.
\end{multline}
Taking the $\inf_{(\alpha,a)\in \Sigma\times\Gamma}$ we obtain the subsolution inequality.

\textbf{Supersolution:}\\
Let $x$ be a minimum for $V+(\phi+g)$ and such that $V+(\phi+g)(x)=0$. As in the subsolution proof we obtain that $x\in D(A^*)$. For $\varepsilon>0$ take $(\alpha_\varepsilon(\cdot), a_\varepsilon(\cdot))$ an $\varepsilon^2$-optimal strategy. We can assume $a(\cdot)$ continuous (it is not hard to see). We call $x(s)$ the trajectory starting from $x$ and subject to $(\alpha_\varepsilon(\cdot),a_\varepsilon(\cdot)$.
Now for $s$ small enough
\bd
V(x(s))+ (\phi+g)(x(s)) \geq V(x) + (\phi+g)(x)
\ed
and thanks to the $\varepsilon^2$-optimality and the Bellman principle we know that
\bd
V(x) + \varepsilon^2 \geq e^{-\rho s} V(x(s)) + \int_0^s e^{-\rho r} L(x(r),\alpha_\varepsilon(r),a_\varepsilon(r)) \ud r
\ed
We take $s=\varepsilon$. Then
\begin{multline}
\frac{1-e^{-\rho \varepsilon}}{\varepsilon} V(x(\varepsilon)) + \frac{\phi(x(\varepsilon)) - \phi(x)}{\varepsilon} +\frac{g(x(\varepsilon)) - g(x)}{\varepsilon}  - \\
- \frac{\int_0^\varepsilon e^{-\rho r} L(x(r),\alpha_\varepsilon(r),a_\varepsilon(r)) \ud r}{\varepsilon} + \frac{\varepsilon^2}{\varepsilon} \geq 0
\end{multline}
in view of Proposition \ref{propg} and Proposition \ref{propphi} we can choose, independently of the control $(\alpha_\varepsilon(\cdot),a_\varepsilon(\cdot))$, a $o(1)$ with $o(1)\xrightarrow{\varepsilon \to 0} 0 \s$ such that:
%fino a qui
\begin{multline}
\rho V(x) + \lla A^* \nabla\phi(x),x\rra + \lla \nabla \phi(x) + \nabla g(x), -\mu x \rra -  \\
- \Bigg ( \frac{\int_0^\varepsilon \lla - \beta \delta_0(\nabla\phi(x) , a_\varepsilon(r) \rra_{\mathbb{R}}+ e^{-\rho r} L(x(r),\alpha_\varepsilon(r),a_\varepsilon(r)) \ud r }{\varepsilon} - \\
- \frac{\int_0^\varepsilon \lla \nabla\phi(x) + \nabla g(x) , \alpha_\varepsilon(r) \rra \ud r }{\varepsilon} \Bigg ) \geq o(1) - \frac{g_0'(|x|)}{|x|}\beta\frac{\|\Gamma\|^2}{2}
\end{multline}
we now take $\inf$ over $a$ and $\alpha$ inside the integral and let $\varepsilon \to 0$ to obtain that
\begin{multline}
\rho V(x) +  \lla A^* \nabla\phi(x),x\rra + \lla \nabla \phi (x) + \nabla g (x), -\mu x \rra-\\
- \inf_{(\alpha,a)\in\Sigma\times\Gamma} \bigg ( 
-\lla \beta\delta_0(\nabla\phi(x)), a \rra_{\mathbb{R}} + L(x,\alpha,a)
- \lla \nabla\phi(x) + \nabla g(x), \alpha \rra
\bigg ) \geq \\
\geq - \frac{g_0'(|x|)}{|x|}\beta\frac{\|\Gamma\|^2}{2}.
\end{multline}
(we observe again that the fact that $o(1)\xrightarrow{\varepsilon \to 0} 0$ uniformly in the control is essential). Therefore $V$ is a solution of the HJB equation.
\end{proof}

%fino a qui -12sera!
\subsection{Uniqueness}
Now we can prove a uniqueness result: we prove the result in the case $\mu\not = 0$. The case $\mu=0$ is simpler and can be proven with small changes in the proof.
\begin{Theorem}
\label{thuniqueness}
Let $L$ satisfy (L1) and (L2) let $\Gamma$ and $\Lambda$ be compact subsets of $\mathbb{R}$. Then given a supersolution $v$ of the HJB equation and a subsolution $u$ we have
\[
u(x)\leq v(x) \;\; for \; every \; x\in\mathcal{H}
\]
In particular there exist at most one solution of the HJB equation
\end{Theorem}
\begin{proof}
We will proceed by contradiction.
Assume that $u$ is a subsolution of the HJB equation and $v$ a supersolution and suppose that there exists $\check x\in\mathcal{H}$ and $\gamma>0$ such that
\bd
(u(\check x)-v(\check x)) > \frac{3\gamma}{\rho} >0
\ed
We take $\gamma <1$. So, taken $\vartheta >0$ small enough we have
\be
\label{duegamma}
u(\check x) - v(\check x) - \vartheta |\check x|^2 > \frac{2 \gamma}{\rho}>0
\ee
We consider $\varepsilon>0$ and $\psi\colon \mathcal{H}\times \mathcal{H} \to \mathbb{R}$ given by
\bd
\psi(x,y) \nd u(x) - v(y) - \frac{1}{2 \varepsilon} |B^{1/2}(x-y)|^2 -\frac{\vartheta}{2} |x|^2 - \frac{\vartheta}{2} |y|^2.
\ed
%By (\ref{duegamma}) we know that 
%\[
%\sup_{(x,y)\in \mathcal{H}\times\mathcal{H}} \psi(x,y)>\frac{2\gamma}{\rho}>0
%\]

Thanks to the boundedness of $u$ and $v$, chosen $\vartheta>0$, there exist $R_\vartheta>0$ such that
\be
\label{(i)}
\psi(0,0) \geq \left ( \sup_{(|x| \geq R_\vartheta)\s or \s (|y| \geq R_\vartheta)} ( \psi(x,y)) \right ) +1
\ee
We set
\bd
S=\left \{ (x,y) \in \mathcal{H} \times \mathcal{H} \s : \s |x|\leq R_\vartheta \s and \s |y|\leq R_\vartheta \right \}
\ed
If we choose $R_\vartheta$ big enough $\check x\in S$. By standard techniques (see \cite{LiYong} page 252) we can find $p$ and $q$ in $\mathcal{H}$ with $|p|<\sigma$ and $|q|<\sigma$ and such that
\bd
(x,y) \mapsto \psi(x,y) - \lla Bp,x \rra - \lla Bq,y \rra
\ed
attains a maximum in $S$. We call $(\bar x, \bar y)$ the point of maximum. If we choose $\sigma$ small enough (for example such that $\sigma \|B\| R_\vartheta <\frac{1}{4}\frac{\gamma}{\rho}$) we know by (\ref{(i)}) that such maximum is in the interior of $S$ and, thanks to (\ref{duegamma}), that
\bd
\psi(\bar x,\bar y) - \lla Bp,\bar x \rra - \lla Bq,\bar y \rra > \frac{3\gamma}{2\rho}.
\ed
Moreover
\be
\label{eqperlafine}
\psi(\bar x,\bar y) > \frac{\gamma}{\rho} \s\s and\s so\s\s u(\bar x) - v(\bar y) >\frac{\gamma}{\rho}.
\ee
We now make some preliminary estimates that we will use in the following:

\textbf{Estimates 1 (on $\varepsilon$)}:\\
We observe that
\[
\left \{
\begin{array}{l}
M\colon [1,0) \to \mathbb{R}\\
M\colon \varepsilon \mapsto \sup_{(x,y)\in \mathcal{H}\times\mathcal{H}} \left ( u(x) - v(y) - \frac{1}{2\varepsilon} \left | B^{1/2} (x-y) \right |^2 \right )
\end{array}
\right .
\]
is non-increasing and bounded and so it admits a limit for $\varepsilon\to 0^+$. So there exists a $\bar\varepsilon>0$ such that, for every $0<\varepsilon_1, \varepsilon_2 \leq \bar\varepsilon$ we have that
\begin{equation}
\label{eqchosenepsilon}
\left | M(\varepsilon_1)- M(\varepsilon_2) \right | <\left ( \frac{\gamma}{16(1+|\mu|)} \right )^2
\end{equation}
We choose now $\varepsilon$, that will be fixed in the sequel of the proof:
\begin{equation}
\label{eqchosenepsilonBIS}
\varepsilon:= min \left \{ \bar\varepsilon, \frac{1}{32C_L^2} \right \}
\end{equation}
($C_L$ is the constant introduced in hypothesys (L1) and (L2)). Now we state and prove a claim that we will use in the following:

\textbf{Claim}\\
\textit{If $\tilde x\in\mathcal{H}$ and $\tilde y\in\mathcal{H}$ satisfy
\begin{equation}
\label{eqhpclaim}
u(\tilde x) - v(\tilde y) - \frac{1}{2\varepsilon} \left | B^{1/2}(\tilde x - \tilde y) \right |^2 \geq M(\varepsilon) - \left ( \frac{\gamma}{16(1+|\mu|)} \right )^2
\end{equation}
then
\begin{equation}
\label{eqstimaBx-By}
\frac{1}{\varepsilon} \left | B^{1/2}(\tilde x - \tilde y) \right |^2 \leq \frac{1}{32} \left ( \frac{\gamma}{(1+|\mu|)} \right )^2
\end{equation}
}
proof of the claim:\\
(We follow the idea used in Lemma 3.2 of \cite{CL6})
\begin{multline}
M(\varepsilon/2) \geq u(\tilde x) - y(\tilde y) - \frac{1}{4\varepsilon} \left | B^{1/2}(\tilde x -\tilde y) \right |^2 =\\
= u(\tilde x) - y(\tilde y) - \frac{1}{2\varepsilon} \left | B^{1/2}(\tilde x -\tilde y) \right |^2 + \frac{1}{4\varepsilon} \left | B^{1/2}(\tilde x -\tilde y) \right |^2 \geq\\
\geq M(\varepsilon) - \left ( \frac{\gamma}{16(1+|\mu|)} \right )^2 + \frac{1}{4\varepsilon} \left | B^{1/2}(\tilde x -\tilde y) \right |^2
\end{multline}
So
\begin{multline}
\frac{1}{4\varepsilon} \left | B^{1/2}(\tilde x -\tilde y) \right |^2 \leq M(\varepsilon/2)- M(\varepsilon) + \left ( \frac{\gamma}{16(1+|\mu|)} \right )^2 \leq \\
\leq \left ( \frac{\gamma}{16(1+|\mu|)} \right )^2 + \left ( \frac{\gamma}{16(1+|\mu|)} \right )^2= 2 \left ( \frac{\gamma}{16(1+|\mu|)} \right )^2
\end{multline}
where the inequality $M(\varepsilon/2)- M(\varepsilon) < \left ( \frac{\gamma}{16(1+|\mu|)} \right )^2$ follows from the definition of $\varepsilon$ (\ref{eqchosenepsilonBIS}) that implies $\varepsilon\leq \bar\varepsilon$ and then the (\ref{eqchosenepsilon}). The claim follows. \qqed

Note that from (\ref{eqchosenepsilonBIS}) we have
\[
\frac{1}{\sqrt{\varepsilon}}\geq 4 \sqrt{2} C_L
\]
and then if $\tilde x$, $\tilde y$ satisfy the hypothesis (\ref{eqhpclaim}) of the claim we have
\begin{equation}
\label{eqstimaCL}
C_L \left |\tilde x -\tilde y \right |_B \leq \frac{\gamma}{32(1+|\mu|)}
\end{equation}

\textbf{Estimates 2 (on $\sigma$)}:\\
We have already imposed $\sigma < \frac{\gamma/\rho}{4 \|B\| R_{\vartheta}}$, we take from now
\begin{equation}
\label{eqsceltasigma}
\sigma= min \left \{ \frac{\gamma}{8 \rho \|B\| R_{\vartheta}}, \vartheta, \frac{\vartheta}{R_{\vartheta}} \right \}
\end{equation}
so that 
\begin{equation}
\label{eqlimsigma}
\sigma \xrightarrow{\vartheta \to 0}0
\end{equation}

and
\begin{equation}
\label{eqlimvarthetaRvartheta}
\sigma R_{\vartheta} \xrightarrow{\vartheta \to 0} 0
\end{equation}
We recall that we have already fixed $\varepsilon$ in (\ref{eqchosenepsilonBIS}). From the choice of $\sigma$ (\ref{eqsceltasigma}) follows that
\begin{equation}
\label{eqstimaBpx}
\left | \lla Bp, \bar x \rra \right | \leq \| B \| \sigma R_{\vartheta} \xrightarrow{\vartheta \to 0}0, \;\;\;\;\; \left | \lla Bq, \bar y \rra \right | \leq \| B \| \sigma R_{\vartheta} \xrightarrow{\vartheta \to 0}0
\end{equation}
Moreover, in view of the continuity of the linear operator $A^*B\colon \mathcal{H}\to\mathcal{H}$ that has norm $\|A^*B\|$, we have
\begin{equation}
\label{eqstimaABpx}
\left | \lla A^*Bp, \bar x \rra \right | \leq \| A^*B \| \sigma R_{\vartheta} \xrightarrow{\vartheta \to 0}0, \;\;\;\;\; \left | \lla A^*Bq, \bar y \rra \right | \leq \| B \| \sigma R_{\vartheta} \xrightarrow{\vartheta \to 0}0
\end{equation}

\textbf{Estimates 3 (on $\vartheta$)}:
One can prove that, fixed $\varepsilon$ we have
\begin{equation}
\label{eqstimaax2}
\vartheta \left | \bar x \right |^2 \xrightarrow{\vartheta \to 0}0, \;\;\;\; \vartheta \left | \bar y \right |^2 \xrightarrow{\vartheta \to 0}0
\end{equation}
(it is a quite standard fact, see for example \cite{CL6}). So
\begin{multline}
\label{eqlimitvartheta}
\lim_{\vartheta \to 0} \left ( \psi(\bar x, \bar y) - \lla Bp, \bar x \rra - \lla Bq, \bar y \rra \right )= \\
=\sup_{(x,y)\in \mathcal{H}\times\mathcal{H}} \left ( x(x)-v(y) - \frac{1}{2\varepsilon} \left | B^{1/2}(x-y) \right |^2 \right )>2\frac{\gamma}{\rho}
\end{multline}
(where the last inequality follows from the (\ref{duegamma})).
In (\ref{eqchosenepsilonBIS}) we fixed $\varepsilon$, in (\ref{eqsceltasigma}) we chose $\sigma$ as function of $\vartheta$. Now we will fix $\vartheta$. 
We begin taking 
\[
\vartheta <\frac{\gamma}{64\beta\|\Gamma\|^2}
\]
so that
\begin{equation}
\label{eqperfinebetathetagamma}
\beta\vartheta\|\Gamma\|^2<\frac{\gamma}{64}
\end{equation}

We know from (\ref{eqstimaBpx}) and (\ref{eqstimaABpx}) that if we choose $\vartheta$ small enough we have
\begin{equation}
\label{eqperfineBpx}
\begin{array}{c}
|\mu|\left | \lla Bp, \bar x \rra \right | <\frac{\gamma}{64}, \;\;\;\;\;\;\;\; |\mu| \left | \lla Bq, \bar y \rra \right | <\frac{\gamma}{64}\\
\\ 
\left | \lla A^*Bp, \bar x \rra \right | <\frac{\gamma}{64}, \;\;\;\;\;\;\;\; \left | \lla A^*Bq, \bar y \rra \right | <\frac{\gamma}{64}
\end{array}
\end{equation}
From (\ref{eqstimaax2}) we know that if we choose $\vartheta$ small enough we have
\begin{equation}
\label{eqperfineax2}
|\mu| \vartheta \left | \bar x \right |^2 <\frac{\gamma}{32}, \;\;\;\; |\mu| \vartheta \left | \bar y \right |^2 <\frac{\gamma}{32}
\end{equation}

Moreover the (\ref{eqstimaax2}) implies also that
\[
\vartheta \left | \bar x \right | \xrightarrow{\vartheta \to 0}0, \;\;\;\; \vartheta \left | \bar y \right | \xrightarrow{\vartheta \to 0}0
\]
and then if we choose $\vartheta$ small enough we have
\begin{equation}
\label{eqperfinesigmathetax}
\vartheta \|\Sigma\| \left ( |\bar x| +|\bar y| \right ) <\frac{\gamma}{32}
\end{equation}

Moreover, in view of (\ref{eqlimitvartheta}) we know that that if we choose $\vartheta$ small enough, $\bar x$ and $\bar y$ satisfy the hypothesis (\ref{eqhpclaim}) of the Claim and then, from the (\ref{eqstimaBx-By}), we have
\begin{equation}
\label{eqperfinestimaBx-By-1}
\frac{1}{\varepsilon} \left | B^{1/2}(\bar x - \bar y) \right |^2 \leq \frac{1}{32} \left ( \frac{\gamma}{(1+|\mu|)} \right )^2 \leq \frac{1}{32} \frac{\gamma}{(1+|\mu|)} \leq \frac{\gamma}{32}
\end{equation}
(where we uses that if $0<a<1$ then $a^2<a$, we recall that we took $0<\gamma<1$). From the (\ref{eqstimaBx-By}) in the same way we obtain
\begin{equation}
\label{eqperfinestimaBx-By-2}
\frac{|\mu|}{\varepsilon} \left | B^{1/2}(\bar x - \bar y) \right |^2 \leq \frac{\gamma}{32}
\end{equation}
and, from (\ref{eqstimaCL}),
\begin{equation}
\label{eqperfinestimaCL}
C_L \left | B^{1/2}(\bar x - \bar y) \right |^2\leq \frac{\gamma}{32}
\end{equation}
Eventually, from (\ref{eqlimsigma}) if we choose $\vartheta$ small enough we have
\begin{equation}
\label{eqperfinestima2Bpsigma}
2 \|B\| |p| \|\Sigma\|\leq 2 \|B\| \sigma \|\Sigma\|\leq \frac{\gamma}{64}
\end{equation}
and
\begin{equation}
\label{eqperfinestima2sigmagamma}
2 \sigma \|\delta_0 \circ B\| \|\Gamma\| \leq \frac{\gamma}{32}
\end{equation}
where we have called $\|\delta_0 \circ B\|$ the norm of the linear continuous functional \mbox{$\delta_0 \circ B \colon \mathcal{H} \to \mathbb{R}$}.

We choose $\vartheta>0$ small enough to satisfy (\ref{eqperfineBpx}), (\ref{eqperfineax2}), (\ref{eqperfinestimaBx-By-1}), (\ref{eqperfinestimaBx-By-2}) (\ref{eqperfinestimaCL}), (\ref{eqperfinestima2Bpsigma}), (\ref{eqperfinestima2sigmagamma}). \qqed

We have finished our preliminary estimates and we come back to the main part of the proof of the theorem. The map
\bd
x \mapsto u(x) - \left ( \frac{1}{2\varepsilon} |B^{1/2}(x-\bar y)|^2 + \frac{\vartheta}{2} |x|^2 + \lla Bp,x\rra \right )
\ed
attains a maximum at $\bar x$ and
\bd
y \mapsto v(y) + \left ( \frac{1}{2\varepsilon} |B^{1/2}(\bar x-y)|^2 + \frac{\vartheta}{2} |y|^2 + \lla Bq,y\rra \right )
\ed
attains a minimum at $\bar y$.

Note that thanks to Proposition \ref{prop0} ${\bar x}$ and ${\bar y}$ are in $D(A^*)$. We can now use the definition of sub- and super-solution (page \pageref{defsubsol}) to obtain that
\begin{multline}
\rho u({\bar x}) - \frac{1}{\varepsilon} \lla A^*B({\bar x}-{\bar y}),{\bar x}\rra - \frac{1}{\varepsilon} \lla B({\bar x}-{\bar y}),-\mu{\bar x}\rra
- \lla A^*B p, {\bar x} \rra - \\
- \lla B p, -\mu{\bar x} \rra - \vartheta\lla \bar x,-\mu \bar x \rra
-\inf_{(\alpha,a)\in\Sigma\times\Gamma} \bigg ( \frac{1}{\varepsilon} \lla \beta \delta_0(B({\bar x}-{\bar y})),a\rra_{\mathbb{R}}
+ \lla \beta \delta_0(Bp),a\rra_{\mathbb{R}} +\\
+ \vartheta \lla \bar x, \alpha \rra + \frac{1}{\varepsilon}\lla B(\bar x - \bar y),\alpha \rra + \lla Bp,\alpha \rra + L({\bar x},\alpha,a) \bigg ) \leq \frac{\vartheta \beta \|\Gamma\|^2}{2} 
\end{multline}
and
\begin{multline}
\rho v({\bar y}) - \frac{1}{\varepsilon} \lla A^*B({\bar x}-{\bar y}),{\bar y}\rra - \frac{1}{\varepsilon} \lla B({\bar x}-{\bar y}),-\mu{\bar y}\rra
+ \lla A^*B q, {\bar y} \rra + \\
+ \lla B q, -\mu{\bar y} \rra + \vartheta\lla \bar y,-\mu \bar y \rra
-\inf_{(\alpha,a)\in\Sigma\times\Gamma} \bigg ( \frac{1}{\varepsilon} \lla \beta \delta_0(B({\bar x}-{\bar y})),a\rra_{\mathbb{R}}
- \lla \beta \delta_0(Bq),a\rra_{\mathbb{R}} -\\
- \vartheta \lla \bar y, \alpha \rra + \frac{1}{\varepsilon}\lla B(\bar x - \bar y),\alpha \rra - \lla Bq,\alpha \rra + L({\bar y},\alpha,a) \bigg ) \geq - \frac{\vartheta \beta \|\Gamma\|^2}{2} 
\end{multline}

Subtracting the above we obtain
\begin{multline}
\label{eqsubtrack}
\rho u (\bar x) - \rho v(\bar y) - \frac{1}{\varepsilon} \lla A^*B(\bar x - \bar y), (\bar x - \bar y) \rra - \\
- \frac{1}{\varepsilon} \lla B(\bar x - \bar y), -\mu (\bar x -\bar y) \rra - \lla A^* B p , \bar x \rra - \lla A^*B q , \bar y \rra - \\
- \lla Bp ,-\mu \bar x \rra - \lla Bq, -\mu \bar y \rra - \vartheta \lla \bar x, -\mu \bar x \rra - \vartheta \lla \bar y , -\mu \bar y \rra-\\
- \inf_{(\alpha,a)\in\Sigma\times\Gamma} \bigg ( \frac{1}{\varepsilon} \lla \beta \delta_0(B({\bar x}-{\bar y})),a\rra_{\mathbb{R}}
+ \lla \beta \delta_0(Bp),a\rra_{\mathbb{R}} +\\
+ \vartheta \lla \bar x, \alpha \rra + \frac{1}{\varepsilon}\lla B(\bar x - \bar y),\alpha \rra + \lla Bp,\alpha \rra + L({\bar x},\alpha,a) \bigg ) +\\
+ \inf_{(\alpha,a)\in\Sigma\times\Gamma} \bigg ( \frac{1}{\varepsilon} \lla \beta \delta_0(B({\bar x}-{\bar y})),a\rra_{\mathbb{R}}
- \lla \beta \delta_0(Bq),a\rra_{\mathbb{R}} -\\
- \vartheta \lla \bar y, \alpha \rra + \frac{1}{\varepsilon}\lla B(\bar x - \bar y),\alpha \rra - \lla Bq,\alpha \rra + L({\bar y},\alpha,a) \bigg )\leq\\
\leq \beta \vartheta \|\Gamma\|^2
\end{multline}
We now note that:\\
\textbf{(A)}: from (\ref{eqRenardy}) $A^*B\leq B$ and then
\bd
- \frac{1}{\varepsilon} \lla A^*B({\bar x}-{\bar y}),({\bar x}-{\bar y})\rra \geq - \frac{1}{\varepsilon} \lla B({\bar x}-{\bar y}),({\bar x}-{\bar y})\rra = - \frac{1}{\varepsilon} |{\bar x}-{\bar y}|^2_B
\ed
\textbf{(B)}: We have
\begin{multline}
- \inf_{(\alpha,a)\in\Sigma\times\Gamma} \bigg ( \frac{1}{\varepsilon} \lla \beta \delta_0(B({\bar x}-{\bar y})),a\rra_{\mathbb{R}}
+ \lla \beta \delta_0(Bp),a\rra_{\mathbb{R}} +\\
+ \vartheta \lla \bar x, \alpha \rra + \frac{1}{\varepsilon}\lla B(\bar x - \bar y),\alpha \rra + \lla Bp,\alpha \rra + L({\bar x},\alpha,a) \bigg ) +\\
+ \inf_{(\alpha,a)\in\Sigma\times\Gamma} \bigg ( \frac{1}{\varepsilon} \lla \beta \delta_0(B({\bar x}-{\bar y})),a\rra_{\mathbb{R}}
- \lla \beta \delta_0(Bq),a\rra_{\mathbb{R}} -\\
- \vartheta \lla \bar y, \alpha \rra + \frac{1}{\varepsilon}\lla B(\bar x - \bar y),\alpha \rra - \lla Bq,\alpha \rra + L({\bar y},\alpha,a) \bigg ) \geq\\
\geq \inf_{(\alpha,a)\in\Sigma\times\Gamma} \bigg ( - \lla \beta \delta_0 (Bp), a \rra_{\mathbb{R}} - \lla \beta \delta_0 (Bq), a \rra_{\mathbb{R}} + L(\bar y,\alpha,a)- L(\bar x,\alpha,a) -\\
- \vartheta \lla \bar y, \alpha \rra - \vartheta \lla \bar x, \alpha \rra - \lla Bq, \alpha \rra - \lla Bp, \alpha \rra \bigg ) \geq\\
\geq \inf_{(\alpha,a)\in\Sigma\times\Gamma} \Big (L(\bar y,\alpha,a)- L(\bar x,\alpha,a) \Big ) - \\
- \sup_{(\alpha,a)\in\Sigma\times\Gamma} \Big (
\lla \beta \delta_0 (Bp), a \rra_{\mathbb{R}} + \lla \beta \delta_0 (Bq), a \rra_{\mathbb{R}} \Big )- \\
- \sup_{(\alpha,a)\in\Sigma\times\Gamma} \Big (\vartheta \lla \bar y, \alpha \rra + \vartheta \lla \bar x, \alpha \rra \Big ) - \sup_{(\alpha,a)\in\Sigma\times\Gamma} \Big (\lla Bq, \alpha \rra + \lla Bp, \alpha \rra \Big ) \geq \\
\geq - C_L |\bar x - \bar y|_B - 2 \sigma \|\delta_0 \circ B\| \|\Gamma\| - \|\Sigma\| \vartheta( |\bar x| +|\bar y|) - 2 \|B\| \sigma \|\Sigma\|
\end{multline}
Thus using (A) and (B) in (\ref{eqsubtrack}) we have
\begin{multline}
\rho \Big ( u({\bar x})-v({\bar y}) \Big ) - \frac{1}{\varepsilon} |{\bar x}-{\bar y}|^2_B - \\
- \frac{\mu}{\varepsilon} \lla B(\bar x - \bar y), - (\bar x -\bar y) \rra - \lla A^* B p , \bar x \rra - \lla A^*B q , \bar y \rra - \\
- \lla Bp ,-\mu \bar x \rra - \lla Bq, -\mu \bar y \rra - \vartheta \lla \bar x, -\mu \bar x \rra - \vartheta \lla \bar y , -\mu \bar y \rra-\\
- C_L |\bar x - \bar y|_B - 2 \sigma \|\delta_0 \circ B\| \|\Gamma\| - \|\Sigma\| \vartheta( |\bar x| +|\bar y|) - 2 \|B\| \sigma \|\Sigma\| -\beta \vartheta \|\Gamma\|^2\leq 0
\end{multline}
using (\ref{eqperfinestimaBx-By-1}), (\ref{eqperfinestimaBx-By-2}), (\ref{eqperfineBpx}), (\ref{eqperfineax2}), (\ref{eqperfinestimaCL}), (\ref{eqperfinestima2sigmagamma}), (\ref{eqperfinesigmathetax}), (\ref{eqperfinestima2Bpsigma}), (\ref{eqperfinebetathetagamma}) we obtain
\begin{multline}
\rho(u({\bar x})-v({\bar y})) - 2\left ( \frac{\gamma}{32} \right )  - 
4\left ( \frac{\gamma}{64} \right ) -
2\left (\frac{\gamma}{32} \right )- 
\frac{\gamma}{32} -
\frac{\gamma}{32} - \frac{\gamma}{32} 
-\frac{\gamma}{64} - \frac{\gamma}{64} \leq 0
\end{multline}
that is
\begin{equation}
\label{equltima}
\rho(u({\bar x})-v({\bar y})) - \frac{1}{2}\gamma \leq 0
\end{equation}
but from the (\ref{eqperlafine}) we have $\rho \big( u({\bar x})-v({\bar y}) \big) > \gamma$ and then we obtain from the (\ref{equltima})
\[
\frac{1}{2}\gamma=\gamma - \frac{1}{2}\gamma <\rho(u({\bar x})-v({\bar y})) - \frac{1}{2}\gamma \leq 0
\]
that is a contradiction because $\gamma>0$ and so the theorem is proven.
\end{proof}

\begin{Remark}
\label{remarkmunondipendedar}
Now we can explain a remark we have done in the introduction: it is difficult to treat with the same arguments the case in which $\mu$ is not a constant but a function of $r$. In the proof of the uniqueness we have to estimate the term $\frac{1}{\varepsilon} \lla B(\bar x - \bar y), -\mu (\bar x -\bar y) \rra$ and we can estimate it because we use the term $\frac{1}{\varepsilon} \left | x-y \right |^2_B$ to penalize the doubling with respect the $B$-norm. 
If we consider the case in which $\mu$ is a function of $r$ such term would appear in the form $\frac{1}{\varepsilon} \lla B(\bar x - \bar y), -\mu(\cdot) (\bar x -\bar y) \rra$ (where $-\mu(\cdot) (\bar x -\bar y)$ is the pointwise product of the $L^{\infty}(0,\bar s)$ function $\mu(\cdot)$ and the $L^2(0,\bar s)$ function $(\bar x -\bar y)$. We do not know how to treat such term.
\end{Remark}

\begin{small}
\bibliographystyle{alpha}
\bibliography{biblioadhoc}
\end{small}

\end{document}